\documentclass[twocolumn]{svjour3}
\usepackage{amssymb,amsmath,color}
\usepackage{subfigure}
\usepackage{graphicx}
\usepackage{warmread}
\usepackage[all,import]{xy}
\smartqed

\newcommand{\norm}[1]{\ensuremath{\left\| #1 \right\|}}

\newcommand{\bracket}[1]{\ensuremath{\left[ #1 \right]}}
\newcommand{\braces}[1]{\ensuremath{\left\{ #1 \right\}}}

\newcommand{\refeqn}[1]{(\ref{eqn:#1})}
\newcommand{\reffig}[1]{Figure \ref{fig:#1}}
\newcommand{\tr}[1]{\mbox{tr}\ensuremath{\negthickspace\bracket{#1}}}
\newcommand{\deriv}[2]{\ensuremath{\frac{\partial #1}{\partial #2}}}
\newcommand{\G}{\ensuremath{\mathsf{G}}}
\newcommand{\SO}{\ensuremath{\mathsf{SO(3)}}}
\newcommand{\T}{\ensuremath{\mathsf{T}}}
\renewcommand{\L}{\ensuremath{\mathsf{L}}}
\newcommand{\so}{\ensuremath{\mathfrak{so}(3)}}

\renewcommand{\Re}{\ensuremath{\mathbb{R}}}

\newcommand{\D}{\ensuremath{\mathbf{D}}}

\newcommand{\Ad}{\ensuremath{\mathrm{Ad}}}

\newcommand{\g}{\ensuremath{\mathfrak{g}}}
\renewcommand{\sb}{\ensuremath{\overline{s}}}
\newcommand{\mub}{\ensuremath{\overline{\mu}}}

\journalname{Nonlinear Dynamics}

\begin{document}

\title{Computational Dynamics of a 3D Elastic String Pendulum Attached to a Rigid Body and an Inertially Fixed Reel Mechanism }

\author{Taeyoung Lee\and
        Melvin Leok\and
        N. Harris McClamroch}

\institute{Taeyoung Lee, Assistant Professor \at
           Department of Mechanical and Aerospace Engineering, Florida Institute of Technology, Melbourne, FL 32901.
           \email{taeyoung@fit.edu}
        \and
           Melvin Leok, Associate Professor,  \at
           Department of Mathematics, University of California, San Diego, CA 92093.
           \email{mleok@ucsd.edu}
        \and
           N. Harris McClamroch, Professor,\at
           Department of Aerospace Engineering, University of
           Michigan, Ann Arbor, MI 48109.
           \email{nhm@engin.umich.edu}
}

\date{September 10, 2009}

\maketitle

\begin{abstract}
A high fidelity model is developed for an elastic string pendulum, one end of which is attached to a rigid body while the other end is attached to an inertially fixed reel mechanism which allows the unstretched length of the string to be dynamically varied. The string is assumed to have distributed mass and elasticity that permits axial deformations. The rigid body is attached to the string at an arbitrary point, and the resulting string pendulum system exhibits nontrivial coupling between the elastic wave propagation in the string and the rigid body dynamics. Variational methods are used to develop coupled ordinary and partial differential equations of motion. Computational methods, referred to as Lie group variational integrators, are then developed, based on a finite element approximation and the use of variational methods in a discrete-time setting to obtain discrete-time equations of motion. This approach preserves the geometry of the configurations, and leads to accurate and efficient algorithms that have guaranteed accuracy properties that make them suitable for many dynamic simulations, especially over long simulation times. Numerical results are presented for typical examples involving a constant length string, string deployment, and string retrieval. These demonstrate the complicated dynamics that arise in a string pendulum from the interaction of the rigid body motion, elastic wave dynamics in the string, and the disturbances  introduced by the reeling mechanism. Such interactions are dynamically important in many engineering problems, but tend be obscured in lower fidelity models.
\keywords{Lagrangian mechanics \and geometric integrator \and variational integrator \and string pendulum \and reel mechanism \and rigid body}
\end{abstract}

\section{Introduction}

\vspace*{0.2cm}

The dynamics of a body connected to a string appear in several engineering problems such as cable cranes, towed underwater vehicles, and tethered spacecraft. Several types of analytical and numerical models have been developed. Lumped mass models, where the string is spatially discretized into connected point masses, are developed in~\cite{ChaOE84,DriLueAOR00,WalPolMC60}. Finite difference methods in both the spatial domain and the time domain are applied in~\cite{KohZhaJEM99,KuhSeiAMC95}.  Finite element discretizations of the weak form of the equations of motion are used in~\cite{KuhSeiAMC95,BucDriTA04}. Variable-length string models also have been developed: a variable length string is modeled based on a continuous plastic impact assumption in~\cite{SteZemAA95,KruPotND06}, and a reel mechanism is considered in~\cite{BanKanJAS82,ManAgrTA05}. But, the reel mechanisms developed in those papers are problematic. In~\cite{BanKanJAS82}, the deployed portion of the string is assumed to move along a fixed line.   The dynamic model of reeling developed in~\cite{ManAgrTA05} is erroneous (this will be discussed further in Section \ref{sec:EL}). Instead of a point mass, a rigid body model is considered in~\cite{KruPotND06}, but this paper does not provide any computational results. Analytical and numerical models of a rigid body connected to an elastic string appear in~\cite{LeeLeoPICDC09}.

The goal of this paper is to develop an analytical model and a numerical algorithm that can be used for simulation of an an elastic string attached to a rigid body and an inertially fixed  reel mechanism, acting under a constant gravitational potential. The string has distributed mass; it can move in a three-dimensional space while deforming axially; the rigid body attached to the string can translate and rotate. We assume that the point where the string is attached to the rigid body is displaced from the center of mass of the rigid body so that there exist nonlinear coupling effects between the string deformation dynamics and the rigid body dynamics. The reel mechanism is an inertially-fixed system consisting of  a cylindrical reel on which the string winds and unwinds and a guide way that acts as a pivot for the deployed portion of the string. The portion of the string on the reel mechanism is assumed to be inextensible. The combined system of the string, the rigid body and the reel mechanism provides a realistic and accurate dynamic model of cable cranes and towing systems. 

In this paper, we first show that the governing equations of motion of the presented string pendulum can be developed according to Hamilton's variational principle. The configuration manifold of the string pendulum is expressed as the product of the real space $\Re$ representing the configuration of the reel mechanism, the space of connected curve segments on $\Re^3$ describing the deployed portion of the string, and the special orthogonal group $\SO$ defining the attitude of the rigid body~\cite{MarRat99}. The variational principle is carefully applied to respect the geometry of the Lie group configuration manifold. We incorporate an additional modification term, referred to as the Carnot energy loss term~\cite{CreJanJAM97}, in the variational principle to take account of the fact that the portion of the string in the reel mechanism is inextensible. The resulting Euler-Lagrange equations are expressed as coupled partial and ordinary differential equations.

The second part of this paper deals with a geometric numerical integrator for the model we presented for the string pendulum. Geometric numerical integration is concerned with developing numerical integrators that preserve geometric features of a system, such as invariants, symmetry, and reversibility~\cite{HaiLub00}. 
The string pendulum is a Lagrangian/Hamiltonian system evolving on a Lie group. When numerically simulating such systems, it is critical to preserve both the symplectic property of Hamiltonian flows and the Lie group structure for numerical accuracy and efficiency~\cite{LeeLeoCMDA07}. A geometric numerical integrator, referred to as a Lie group variational integrator, has been developed for a Hamiltonian system on an arbitrary Lie group and it has been applied to several multibody systems ranging from binary asteroids to articulated rigid bodies and magnetic systems in~\cite{LeeLeoCMAME07,Lee08}. 

This paper develops a Lie group variational integrator for the proposed string pendulum model. This extends the results presented in~\cite{Lee08} by incorporating deformation of the string using a finite element model and by including a discrete-time Carnot energy loss term. The proposed geometric numerical integrator preserves the symplectic structure and momentum maps, and exhibits desirable energy conservation properties. It also respects the Lie group structure of the configuration manifold, and avoids the singularities and computational complexities associated with the use of local coordinates, explicit constraints or projection. As a result, this computational approach can represent arbitrary translations and rotations of the rigid body and large deformations of the string.

In summary, this paper develops an analytical model and a geometric numerical integrator for a string pendulum attached to a rigid body and a reel mechanism. These provide a realistic mathematical model for tethered systems and a reliable numerical simulation tool that characterizes the nonlinear coupling between the string dynamics, the rigid body dynamics, and the reel mechanism accurately. The proposed high-fidelity computational framework can be naturally extended to formulating and solving control problems associated with string deployment, retrieval, and vibration suppression as in~\cite{LeLeMc2008a}. 

This paper is organized as follows. A string pendulum is described and the corresponding Euler-Lagrange equations are presented in Section~\ref{sec:EL}. A Lie group variational integrator is derived in Section~\ref{sec:LGVI}, followed by numerical examples and conclusions in Section~\ref{sec:NE} and \ref{sec:CON}.

\section{Euler-Lagrange Equations}\label{sec:EL}

\subsection{String Pendulum Model}

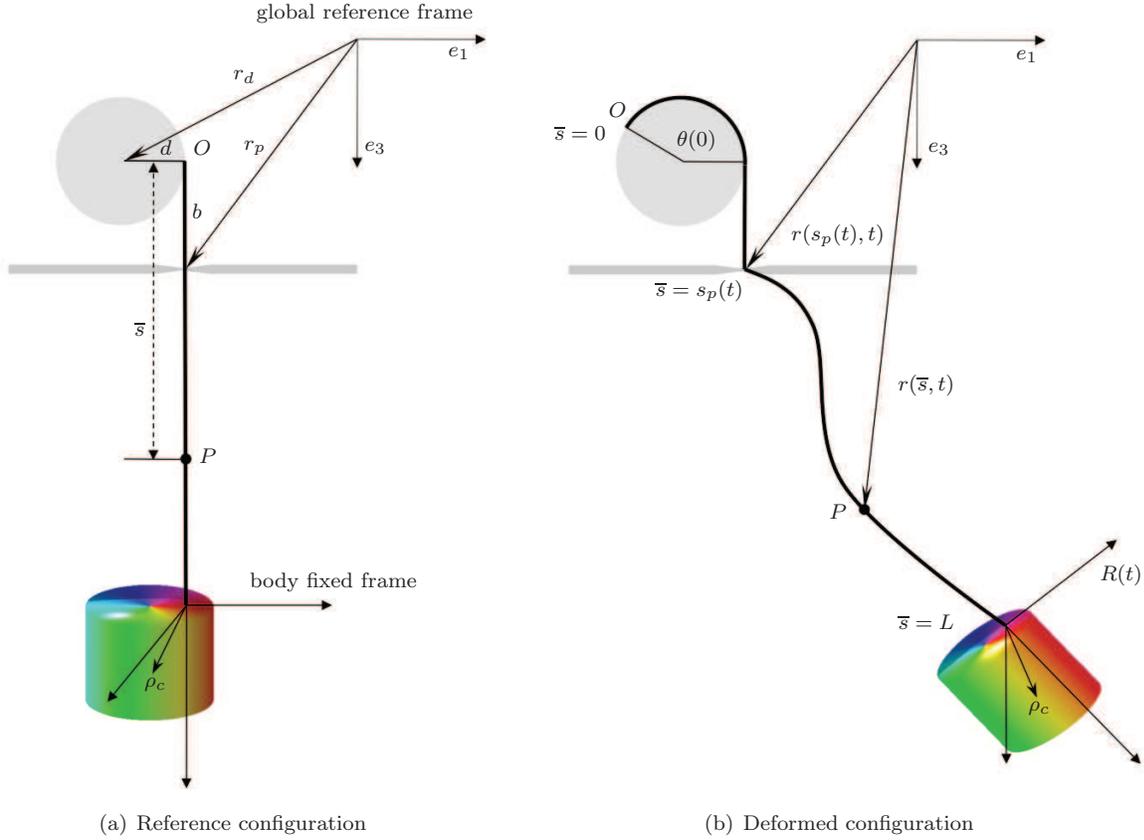
\begin{figure*}
\centerline{
\subfigure[Reference configuration]{
\renewcommand{\xyWARMinclude}[1]{\includegraphics[width=0.40\textwidth]{#1}}
{\footnotesize\selectfont
$$\begin{xy}
\xyWARMprocessEPS{Lee1a}{eps}
\xyMarkedImport{}
\xyMarkedMathPoints{1-15}
\end{xy}\vspace*{-0.35cm}
$$}}
\subfigure[Deformed configuration]{
{\footnotesize\selectfont
\renewcommand{\xyWARMinclude}[1]{\includegraphics[width=0.481\textwidth]{#1}}
$$\begin{xy}
\xyWARMprocessEPS{Lee1b}{eps}
\xyMarkedImport{}
\xyMarkedMathPoints{1-15}
\end{xy}\vspace*{-0.35cm}
$$}}}
\caption{String Pendulum Model}\label{fig:SP}
\end{figure*}

Consider a string that is composed of mass elements distributed along a curve. The string mass elements can translate in a three-dimensional space, and it is deformable along its axial direction. The bending stiffness of the string is not considered as the diameter of the string is assumed to be negligible compared to its length. The free end of the string is attached to a rigid body that can translate and rotate, and the point where the string is attached to the rigid body is displaced from the center of mass of the rigid body so that the dynamics of the rigid body is coupled to the string deformations and displacements. The other end of the string is connected to an inertially-fixed reel mechanism composed of a drum and a guide way. The string is wound around the drum at a constant radius, and the string on the drum and in the guide way is assumed to be inextensible. A control moment is applied at the rotating drum. This system of the string, the rigid body, and the reel mechanism, acting under a constant gravitational potential, is referred to as a string pendulum. This is illustrated in \reffig{SP}.

We choose an inertially fixed  reference frame and a body-fixed frame. The origin of the body-fixed frame is located at the end of the string where the string is attached to the rigid body, and it is fixed to the rigid body. Since the string is extensible, we need to distinguish between the arc length for the stretched deformed configuration and the arc length for the unstretched reference configuration. Define

\begin{list}{}
{\setlength{\leftmargin}{3.0cm}\setlength{\itemindent}{-2.6cm}
\setlength{\parsep}{0cm}\setlength{\itemsep}{0.0cm}\setlength{\parskip}{0cm}}
\item\makebox[2.5cm][l]{$r_d\in\Re^3$} the location of the origin of the axis of the  drum
\item\makebox[2.5cm][l]{$d\in\Re$} the radius of the drum
\item\makebox[2.5cm][l]{$I_d=\kappa_d d^2\in\Re$} the rotational  inertia of the drum for $\kappa_d\in\Re$
\item\makebox[2.5cm][l]{$b\in\Re$} the length from the drum to the guide way
\item\makebox[2.5cm][l]{$u\in\Re$} the control moment applied at drum
\item\makebox[2.5cm][l]{$L\in\Re$} the total unstretched length of the string
\item\makebox[2.5cm][l]{$\overline\mu\in\Re$} the mass of the string per the unit unstretched length
\item\makebox[2.5cm][l]{$O$} the point at which the string is attached to the drum
\item\makebox[2.5cm][l]{$\overline s\in[0,L]$} the unstretched arc length of the string between the point $O$ and a material point $P$ on the string
\item\makebox[2.5cm][l]{$s(\overline s,t)\in\Re$} the stretched arc length to a material point $P$
\item\makebox[2.5cm][l]{$s_p(t)\in[b,L]$} the arc length of the string between the point $O$ and the material point on the string located at the guide way entrance
\item\makebox[2.5cm][l]{$r(\overline s,t)\in\Re^3$} the deformed location of a material point $P$
\item\makebox[2.5cm][l]{$\theta(\sb)\in\Re$} $\theta=((s_p-b)-\sb)/d$ for $\sb\in[0,s_p-b]$
\item\makebox[2.5cm][l]{$M\in\Re$} the mass of the rigid body
\item\makebox[2.5cm][l]{$J\in\Re^3$} the inertia matrix of the rigid body with respect to the body fixed frame
\item\makebox[2.5cm][l]{$\rho_c\in\Re^3$} the vector from the origin of the body fixed frame to the center of mass of the rigid body represented in the body fixed frame
\item\makebox[2.5cm][l]{$R\in\SO$} the rotation matrix from the body fixed frame to the reference frame
\item\makebox[2.5cm][l]{$\Omega\in\Re^3$} the angular velocity of the rigid body represented in the body fixed frame
\end{list}

The configuration of the string on the drum and in the guide way is completely determined by the variable $s_p(t)$, since the string there is inextensible. The configurations of the deployed portion of the string and the rigid body are described by the curve $r(\sb,t)$ for $\sb\in[s_p,L]$, and the rotation matrix $R\in\SO$, respectively, where the special orthogonal group is $\SO=\{R\in\Re^{3\times 3}\,|\, R^TR=I,\mathrm{det}[R]=1\}$. Therefore, the configuration manifold of the string pendulum is the product of the real space $\Re$, the space of connected curves on $\Re^3$, and the special orthogonal group $\SO$.

The attitude kinematics equation of the rigid body is given by
\begin{align}
\dot R = R\hat \Omega,\label{eqn:AK}
\end{align}
where the \textit{hat map} $\hat\cdot : \Re^3\rightarrow\so$ is defined by the condition that $\hat x y = x\times y$ for any $x,y\in\Re^3$. Since $\hat x$ is a $3\times 3$ skew-symmetric matrix, we have $\hat x^T=-\hat x$. The inverse map of the hat map is referred to as the \textit{vee map}: $(\cdot)^\vee:\so\rightarrow\Re^3$.  

\subsection{Lagrangian}

We develop Euler-Lagrange equations for the string pendulum according to Hamilton's variational principle. The Lagrangian of the string pendulum is derived, and the corresponding action integral is defined. Due to the unique dynamic characteristics of the string pendulum, the variation of the action integral should be carefully developed: (i) since the unstretched length of the deployed portion of the string is not fixed, when deriving the variation of the corresponding part of the action integral, we need to apply Green's theorem; (ii) since the attitude of the rigid body is represented in the special orthogonal group, the variation of rotation matrices are carefully expressed by using the exponential map~\cite{Lee08,LeeLeoCMAME07}; (iii) since the portion of the string on the guide way and the drum is inextensible, the velocity of the string is not continuous at the guide way entrance. To take account of the effect of this velocity discontinuity, an additional modification term, referred to as a Carnot energy loss term is incorporated~\cite{CreJanJAM97}. Then, Euler-Lagrange equations are derived according to Hamilton's principle, and they are expressed as coupled ordinary and partial differential equations.

\paragraph{Lagrangian}

The total kinetic energy is composed of the kinetic energy of the portion of the string on the drum and the guide way $T_r$, the kinetic energy of the deployed portion of the string $T_s$, and the kinetic energy of the rigid body $T_b$. The kinetic energy $T_r$ can be written as
\begin{align*}
T_{r} = \int_{0}^{s_p} \frac{1}{2} \mub \dot r(\sb)\cdot \dot r(\sb)\, d\sb
+\frac{1}{2}I_d \dot \theta(0)^2,
\end{align*}
where the dot represents the partial derivative with respect to time. Here, the dependency of variables on time $t$ is omitted for simplicity, i.e. $r(\sb)=r(\sb,t)$. The velocity of the string in the reel mechanism is equal to $\dot s_p$ as the string is inextensible. From the definitions, we have $\dot\theta(0)=\dot s_p / d$, and $I_d=\kappa_d d^2$. Then, the kinetic energy $T_r$ can be written as
\begin{align}
T_r = \frac{1}{2} (\mub s_p + \kappa_d) \dot s_p^2.\label{eqn:Tr}
\end{align}
The kinetic energy of the deployed portion of the string is given by
\begin{align}
T_{s} = \int_{s_p}^{L} \frac{1}{2} \mub \dot r(\sb)\cdot \dot r(\sb)\, d\sb.\label{eqn:Ts}
\end{align}
Let $\rho\in\Re^3$ be the vector from the free end of the string $r(L)$ to a mass element of the rigid body,  expressed in the body fixed frame. The location of the mass element in the reference  frame is given by $r(L)+R\rho$. Then, the kinetic energy of the rigid body is given by
\begin{align}
T_{b} &= \int_{\text{body}} \frac{1}{2} \|\dot r(L) + R\hat\Omega \rho\|^2\,dm\nonumber\\
&=\frac{1}{2} M \dot r(L)\cdot \dot r(L) + M\dot r(L) \cdot R\hat\Omega \rho_c + \frac{1}{2}\Omega \cdot J\Omega,\label{eqn:Tb}
\end{align}
where $J=\int -\hat\rho^2\,dm$ is the inertia matrix of the rigid body in the body fixed frame. 

Now we obtain expressions for the potential energy of each part. The gravitational potential energy of the portion of the string on the drum and the guide way is given by
\begin{align}
V_r & = -\int_{0}^{s_p-b} \mub g (r_d\cdot e_3 - d\sin\theta ) \,d\sb\nonumber\\
& =-\mub g \braces{ (s_p-d)\, r_d\cdot e_3 +d^2 (\cos ((s_p-b)/d)-1) }.\label{eqn:Vr}
\end{align}
The strain of the string at a material point located at $r(\sb)$ is given by
\begin{align*}
\epsilon = \lim_{\Delta \sb\rightarrow 0} \frac{\Delta s(\sb)-\Delta \sb}{\Delta \sb}
= s'(\sb) -1,
\end{align*}
where $(\;)'$ denote the partial derivative with respect to $\sb$. The tangent vector at the material point is given by 
\begin{align*}
e_t= \deriv{r(\sb)}{s} = \deriv{r(\sb)}{\sb}\deriv{\sb}{s(\sb)} = \frac{r'(\sb)}{s'(\sb)}.
\end{align*}
Since this tangent vector has unit length, we have $s'(\sb)=\norm{r'(\sb)}$. Therefore, the strain of the string is given by $\epsilon = \norm{r'(\sb)}-1$. The potential energy of the deployed portion of the string is composed of the elastic potential and the gravitational potential energy:
\begin{align}
V_{s} = \int^{L}_{s_p} \frac{1}{2}EA (\norm{r'(\sb)}-1)^2 -\mub g r(\sb)\cdot e_3 \,d\sb,\label{eqn:Vs}
\end{align}
where $E$ and $A$ denote the Young's modulus and the cross sectional area of the string, respectively. The gravitational potential of the rigid body is given by
\begin{align}
V_b & = - M g (r(L)+R\rho_c)\cdot e_3. \label{eqn:Vb}
\end{align}

In summary, the Lagrangian of the string pendulum is given by
\begin{align}
L = (T_r-V_r) + (T_s-V_s) + (T_b-V_b) = L_r + L_s + L_b.
\end{align}

\subsection{Variational Approach}

\paragraph{Action Integral} The action integral is defined by
\begin{align}
\mathfrak{G} = \int_{t_0}^{t_f} L_r+L_s+L_b\,dt = \mathfrak{G}_r + \mathfrak{G}_s+\mathfrak{G}_b.\label{eqn:G}
\end{align}
We find expressions for the variation of each term of the action integral. 

\paragraph{Variation of $\mathfrak{G}_r$}

From \refeqn{Tr} and \refeqn{Vr}, the variation of $\mathfrak{G}_r$ is given by
\begin{align}
\delta \mathfrak{G}_r & = \int_{t_0}^{t_f} \bigg\{-(\mub s_p + \kappa_d)\ddot s_p -\frac{1}{2}\mub \dot s_p^2 + \mub g \,(r_d\cdot e_3)\nonumber\\
& \quad - \mub g d\sin((s_p-b)/d)\bigg\}\,\delta s_p \,dt, \label{eqn:delGr}
\end{align}
where we used integration by parts.

\paragraph{Variation of $\mathfrak{G}_s$}
From \refeqn{Ts}, \refeqn{Vs}, \refeqn{G}, the second term of the action integral $\mathfrak{G}_s$ is a double integral on $(t,\sb)\in[t_0,t_f]\times[s_p(t),L]$. Since the variable $s_p(t)$ is dependent on the time $t$, the variation of $\mathfrak{G}_s$ should take into account the variation of $s_p(t)$:
\begin{align}
\delta& \mathfrak{G}_s = \int_{t_0}^{t_f} \int_{s_p(t)}^L 
\mub \dot r(\sb)\cdot \delta\dot r(\sb) \nonumber\\
&\quad -EA \frac{\norm{r'(\sb)}-1}{\norm{r'(\sb)}}r'(\sb)\cdot \delta r'(\sb) 
+\mub g e_3 \cdot \delta r(s)\,d\sb\,dt\nonumber\\
&\quad - \int_{t_0}^{t_f} \bigg\{\frac{1}{2}\mub\dot r(s_p^+) \cdot \dot r(s_p^+)
-\frac{1}{2}EA(\norm{r'(s_p^+)}-1)^2\nonumber\\
&\quad+\mub g r(s_p)\cdot e_3\bigg\}\,\delta s_p\,dt,\label{eqn:delGs0}
\end{align}
where $r(s_p^+)$ represents the material point of the string located just outside the guide way. 

Now we focus on the first term of \refeqn{delGs0}. Here, we cannot apply integration by parts at time $t$, since the order of the integrals in \refeqn{delGs0} cannot be interchanged due to the time dependence in the variable $s_p(t)$.  Instead, we use Green's theorem,
\begin{align}
\oint_\mathcal{B} \dot r(\sb)\cdot\delta r(\sb) \,d\sb 
= \int_{t_0}^{t_f} \int_{s_p(t)}^L \frac{d}{dt}(\dot r(\sb)\cdot\delta r(\sb))\,d\sb dt,\label{eqn:Green1}
\end{align}
where $\oint_\mathcal{B}$ represents the counterclockwise line integral on the boundary $\mathcal{B}$ of the region $[t_0,t_f]\times[s_p(t),L]$. The boundary $\mathcal{B}$ is composed of four lines: $(t=t_0,\sb\in[s_p(t_0),L])$, $(t=t_f,\sb\in[s_p(t_f),L])$, $(t\in[t_0,t_f],\sb=s_p(t))$, and $(t\in[t_0,t_f],\sb=L)$. For the first two lines, $\delta r(\sb)=0$ since $t=t_0,t_f$. For the last line, $d\sb=0$ since $\sb$ is fixed. Thus, parameterizing the third line by $t$, we obtain 
\begin{align*}
\oint_\mathcal{B} \dot r(\sb)\cdot\delta r(\sb) \,d\sb 
=\int_{t_0}^{t_f} \dot r(s_p(t))\cdot \delta r(s_p(t))\, \dot s_p(t)\, dt.
\end{align*}
Substituting this into \refeqn{Green1}
and rearranging, the first term of \refeqn{delGs0} is given by
\begin{align}
&\int_{t_0}^{t_f} \int_{s_p}^L \dot r(\sb)\cdot \delta \dot r(\sb)\nonumber\\& =
\int_{t_0}^{t_f} \bracket{\int_{s_p}^L -\ddot r(\sb)\cdot \delta r(\sb) \,d\sb 
+ \dot r(s_p)\cdot \delta r(s_p)\, \dot s_p\,} dt.\label{eqn:delGs01}
\end{align}
Substituting this into \refeqn{delGs0}, and using integration by parts with respect to $\sb$ for the second term of \refeqn{delGs0}, the variation of $\mathfrak{G}_s$ can be written as
\begin{align*}
\delta& \mathfrak{G}_s  = \int_{t_0}^{t_f} \int_{s_p}^L 
\braces{-\mub\ddot r(\sb)+F'(\sb) +\mub g e_3} \cdot \delta r(\sb)\,d\sb\,dt\nonumber\\
&+ \int_{t_0}^{t_f} \bigg\{-\frac{1}{2}\mub\dot r(s_p^+) \cdot \dot r(s_p^+)
+\frac{1}{2}EA(\norm{r'(s_p^+)}-1)^2\nonumber\\
& \quad -\mub g r(s_p)\cdot e_3\bigg\}\,\delta s_p + 
 \mub\dot r(s_p^+)\cdot \delta r(s_p^+)\dot s_p\nonumber\\
& \quad-F(L)\cdot \delta r(L) + F(s_p)\cdot \delta r(s_p^+)\,dt,
\end{align*}
where $F(\sb) = EA\frac{\norm{r'(\sb)}-1}{\norm{r'(\sb)}}r'(\sb)$ represents the tension of the string.

We simplify this using the boundary  condition at the guide way. The location of the guide way entrance is given by $r_p=r(s_p(t),t)$. Since the location is inertially fixed, we have $\delta r_p = \delta r(s_p^+) + r'(s_p^+)\delta s_p=0$, and $\dot r_p = \dot r(s_p^+) + r'(s_p^+)\dot s_p=0$. Substituting these, we obtain 
{\allowdisplaybreaks
\begin{align}
\delta & \mathfrak{G}_s = \int_{t_0}^{t_f} \int_{s_p(t)}^L 
\braces{-\mub\ddot r(\sb)+F'(\sb) +\mub g e_3} \cdot \delta r(\sb)\,d\sb\,dt\nonumber\\
&\quad  + \int_{t_0}^{t_f} \bigg\{\frac{1}{2}\mub\norm{r'(s_p^+)}^2 \dot s_p^2
+\frac{1}{2}EA(\norm{r'(s_p^+)}-1)^2\nonumber\\
&\quad -\mub g r(s_p)\cdot e_3\bigg\}\,\delta s_p\,dt\nonumber\\
&\quad -\int_{t_0}^{t_f}F(s_p)\cdot r'(s_p^+)\,\delta s_p + F(L)\cdot \delta r(L) \,dt.\label{eqn:delGs}
\end{align}}

\paragraph{Variation of $\mathfrak{G}_b$} 
From \refeqn{Tb}, \refeqn{Vb}, the variation of $\mathfrak{G}_b$ is given by
\begin{align}
\delta & \mathfrak{G}_b = \int_{t_0}^{t_f} \{M \dot r(L)+MR\hat\Omega\rho_c\}\cdot \delta\dot r(L)\nonumber\\
&\quad+\{J\Omega+M\hat\rho_c R^T \dot r(L)\}\cdot\delta\Omega
+ M\dot r(L)\cdot \delta R\hat\Omega \rho_c \nonumber\\
&\quad+ M g e_3\cdot \delta r(L) + M g e_3 \cdot \delta R\rho_c\,dt.\label{eqn:delGb0}
\end{align}

The attitude of the rigid body is represented by the rotation matrix $R\in\SO$. Therefore, the variation of the rotation matrix should be consistent with the geometry of the special orthogonal group. In~\cite{Lee08,LeeLeoCMAME07}, it is expressed in terms of the exponential map as
\begin{align}
\delta R =\frac{d}{d\epsilon}\bigg|_{\epsilon=0} R^\epsilon= \frac{d}{d\epsilon}\bigg|_{\epsilon=0} R \exp \epsilon\hat\eta = R \hat\eta \label{eqn:delR}
\end{align}
for $\eta\in\Re^3$. The key idea is expressing the variation of a Lie group element in terms of a Lie algebra element. This is desirable since the Lie algebra $\so$ of the special orthogonal group, represented by $3\times 3$ skew symmetric matrices, is isomorphic as a Lie algebra to $\Re^3$ using the vee map. As a result, the variation of the three-dimenstional rotation matrix $R$ is expressed in terms of a vector $\eta\in\Re^3$. We can directly show that \refeqn{delR} satisfies $\delta(R^TR)=\delta R^T R + R^T \delta R =-\hat\eta+\hat\eta=0$. The corresponding variation of the angular velocity is obtained from the kinematics equation \refeqn{AK}:
\begin{align}
\delta\hat\Omega = \frac{d}{d\epsilon}\bigg|_{\epsilon=0} (R^\epsilon)^T \dot R^\epsilon = (\dot\eta + \Omega\times \eta)^\wedge.\label{eqn:delw}
\end{align}

Substituting \refeqn{delR}, \refeqn{delw} into \refeqn{delGb0}, and using integration by parts for $\eta$, we obtain
\begin{align}
\delta & \mathfrak{G}_b  = \int_{t_0}^{t_f} 
M\{-\ddot r(L)-R\hat\Omega^2\rho_c-R\hat{\dot\Omega}\rho_c +ge_3\}\cdot  r(L)\nonumber\\
& +\{-J\dot\Omega-M\hat\rho_c R^T\ddot r(L)
-\hat\Omega J\Omega+Mg\hat\rho_c R^T e_3\}\cdot\eta\,dt.\label{eqn:delGb}
\end{align}

\paragraph{Variation of $\mathfrak{G}$}

From \refeqn{delGr}, \refeqn{delGs}, \refeqn{delGb}, the variation of the action integral is given by
\begin{align}
\delta\mathfrak{G} = \delta\mathfrak{G}_r + \delta\mathfrak{G}_s + \delta\mathfrak{G}_b.\label{eqn:delG}
\end{align}

\paragraph{Variational Principle with Discontinuity}

Let $r_p=r(s_p(t),t)$ be the location of the pivot in the reference frame. Since it is fixed, we have $\dot r_p = \dot r(s_p,t)+r'(s_p,t)\dot s_p=0$. Let $r(s_p^-)$, and $r(s_p^+)$ be the material point of the string just inside the guide way, and the material point just outside the guide way, respectively. Since the string is inextensible inside the guide way, $\|r'(s_p^-)\|=1$. Since the string is extensible outside the guide way, $\|r'(s_p^+)\|=1+\epsilon^+$, where $\epsilon^+$ represents the strain of the string just outside the guide way. Using these, the speeds of the string at those points are given by
\begin{align*}
\|\dot r(s_p^-)\| & = \|-r'(s_p^-)\dot s_p\| = |\dot s_p|\\
\|\dot r(s_p^+)\| & = \|-r'(s_p^+)\dot s_p\| = (1+\epsilon^+)|\dot s_p|.
\end{align*}
Therefore, the speed of the string changes instantaneously by the amount $\epsilon^+|\dot s_p|$ at the guide way.

Due to this velocity and strain discontinuity, the variation of the action integral is not equal to the negative of the virtual work done by the external control moment $u$ at the drum. In order to derive equations of motion using Hamilton's principle, an additional term $Q$, referred to as Carnot energy loss term should be introduced \cite{CreJanJAM97,KruPotND06}. The resulting variational principle is given by
\begin{align}
\delta\mathfrak G + \int_{t_0}^{t_f} ( Q + u/d )\delta s_p\, dt =0.\label{eqn:VP}
\end{align}
The corresponding time rate of change of the total energy is given by $\dot E = (Q + u/d ) \dot s_p$, where the first term $Q\dot s_p$ represents the energy dissipation rate due to the velocity and strain discontinuity.

Consider the infinitesimal mass element $dm = \mub \dot s_p dt$ located just outside the guide way. Without loss of generality, we assume that $\dot s_p > 0$ (retrieval case). The motion of this mass element moving with the velocity $(1+\epsilon^+)\dot s_p$ can be considered as a plastic impact  into the the portion of the string on the guide way moving with velocity $\dot s_p$. The corresponding energy dissipation rate is given by
\begin{align}
Q\dot s_p = -\frac{1}{2}\mub (\epsilon^+)^2 \dot s_p^3 - \frac{1}{2} EA (\epsilon^+)^2 \dot s_p.\label{eqn:dotE}
\end{align}
(See \cite{CreJanJAM97,KruPotND06}). Dividing both side by $\dot s_p$, we obtain the expression for the Carnot energy loss term $Q$.

\subsection{Euler-Lagrange Equations} 

Substituting \refeqn{delG}, \refeqn{dotE} into \refeqn{VP}, we obtain Euler-Lagrange equations for the string pendulum:
{\allowdisplaybreaks
\begin{gather}
\begin{aligned}
-(\mub s_p & + \kappa_d)\ddot s_p  + \mub g \,(r_d-r_p)\cdot e_3  - \mub g d\sin((s_p-b)/d)\\&-F(s_p^+)\cdot r'(s_p^+)+\mub(\norm{r'(s_p^+)}-1) \dot s_p^2+\frac{u}{d}=0,\end{aligned}\label{eqn:EL0}\\
-\mub\ddot r(\sb)+F'(\sb) +\mub g e_3=0,\quad \sb\in[s_p,L],\label{eqn:EL1}\\
-M\ddot r(L)-MR\hat\Omega^2\rho_c-MR\hat{\dot\Omega}\rho_c +Mge_3-F(L)=0,\label{eqn:EL2}\\
-J\dot\Omega-M\hat\rho_c R^T\ddot r(L)
-\hat\Omega J\Omega+Mg\hat\rho_c R^T e_3=0,\label{eqn:EL3}
\end{gather}}
where $F(\sb)=EA\frac{\norm{r'(\sb)}-1}{\norm{r'(\sb)}}r'(\sb)$ is the tension of the string. These are coupled ordinary and partial differential equations. The motion of the reel mechanism and the deployed portion of the string are described by \refeqn{EL0} and \refeqn{EL1}, respectively. The translational and rotational dynamics of the rigid body are determined by \refeqn{EL2}, \refeqn{EL3}. All of these equations are coupled. 

In \refeqn{EL0}, the fifth term, $\mub\epsilon^+ \dot s_p^2$, represents the effect of the velocity discontinuity. Note that this term vanishes if the deployed portion of the string is also inextensible, i.e. $\norm{r'(\sb)}=1$ or $\epsilon^+=0$. A similar expression is developed in \cite{ManAgrTA05} from momentum balance, but their expression  is erroneous. 

\paragraph{Special Cases}

Suppose that the length of the string on the reel mechanism is fixed, i.e. $s_p(t)\equiv s_p(t_0)$ for any $t>t_0$. Then, the equations of motion \refeqn{EL1}-\refeqn{EL3} describe the dynamics of an elastic  string pendulum model with a fixed unstretched length, which is studied in~\cite{LeeLeoPICDC09}. In this case, the total energy and the total angular momentum about the gravity direction $e_3$ are conserved:
\begin{align*}
E & = (T_r+V_r) + (T_s + V_s) + (T_b + V_b),\\
\pi_3 & = e_3 \cdot \bigg[ \int_{s_p}^L \mub r(\sb)\times \dot r(\sb)\,d\sb\\
& +Mr(L)\times(\dot r(L)+R\hat\Omega\rho_c)-M\dot r(L)\times R\rho_c + J\Omega\bigg].
\end{align*}

If we choose $\rho_c=0$, then the rotational dynamics of the rigid body \refeqn{EL3} is decoupled from the other equations. In this case, \refeqn{EL0}-\refeqn{EL2} describe the dynamics of an elastic string attached to a point mass and a reel mechanism. 

\section{Lie Group Variational Integrator}\label{sec:LGVI}

Geometric numerical integration deals with numerical integration methods that preserve geometric properties of a dynamic system, such as invariants, symmetries, reversibility, or structure of the configuration manifold~\cite{HaiLub00,LeiRei04}. The geometric structure of a dynamic system determines its qualitative dynamical behavior, and therefore, the geometric structure-preserving properties of a geometric numerical integrator play an important role in the qualitatively accurate computation of long-term dynamics. The continuous-time Euler-Lagrange equations developed in the previous section provide an analytical model for a string pendulum. However, the popular finite difference approximations or finite element approximations of those equations using a general purpose numerical integrator may not accurately preserve the geometric properties of the system~\cite{HaiLub00}.

Variational integrators provide a systematic method of developing geometric numerical integrators for Lagrangian/Hamiltonian systems~\cite{MarWesAN01}. Discrete-time Euler-Lagrange equations, referred to as variational integrators, are constructed by discretizing Hamilton's principle, rather than discretizing the continuous-time Euler-Lagrange equations using finite difference approximations. This is in contrast to the conventional viewpoint that a numerical integrator of a dynamic system is a discrete approximation of its continuous-time equations of motion. As it is derived from a discrete analogue of Hamilton's principle, it preserves symplecticity and the momentum map, and it exhibits good total energy behavior for an extremely long time period. 

On the other hand, Lie group methods conserve the structure of a Lie group configuration manifold as it updates a group element using the group operation~\cite{IseMunAN00}. As opposed to computational methods based on local coordinates, projections, or constraints, this approach preserves the group structure naturally without any singularities associated with local coordinates or the additional computational overhead introduced by constraints.

These two methods have been unified to obtain a Lie group variational integrator for Lagrangian/Hamiltonian systems evolving on a Lie group~\cite{Lee08}. This geometric integrator preserves symplecticity and group structure of those systems concurrently. It has been shown that this property is critical for accurate and efficient simulations of rigid body dynamics~\cite{LeeLeoCMDA07}. This is particularly useful for dynamic simulation of a string pendulum that undergoes large displacements, deformation, and rotations over an exponentially long time period.

In this section, we develop a Lie group variational integrator for a string pendulum. We first construct a discretized string pendulum model, and derive an expression for a discrete Lagrangian, which is  substituted into discrete-time Euler-Lagrange equations on a Lie group.

\subsection{Discretized String Pendulum Model}

Let $h>0$ be a fixed time step. The value of variables at $t=t_0+kh$ is denoted by a subscript $k$.  We discretize the deployed portion of the string using $N$ identical line elements. Since the unstretched length of the deployed portion of the string is $L-s_{p_k}$, the unstretched length of each element is $l_{k} = \frac{L-s_{p_k}}{N}$. Let the subscript $a$ denote the variables related to the $a$-th element. The natural coordinate of the $a$-th element is defined by
\begin{align}
\zeta_{k,a} (\sb) = \frac{(\sb-s_{p_k})-(a-1)l_{k}}{l_{k}}
\end{align}
for $\sb\in[s_{p_k}+(a-1)l_{k},s_{p_k}+al_{k}]$. This varies between 0 and 1 for the $a$-th element. Let $S_0,S_1$ be shape functions given by $S_0(\zeta)=1-\zeta$, and $S_1(\zeta)=\zeta$. These shape functions are also referred to as  \textit{tent functions}. Define $q_{k,a}$ to be the relative location of a string element with respect to the guide way entrance, i.e $q_{k,a}=r_{k,a}-r_p$. 

Using this finite element model, the position vector $r(\sb,t)$ of a material point in the $a$-th element is approximated as follows:
\begin{align}
r_k(\sb)= S_0(\zeta_{k,a}) q_{k,a} + S_1(\zeta_{k,a}) q_{k,a+1}+r_p.\label{eqn:rksb}
\end{align}
Therefore, a configuration of the presented discretized string pendulum at $t=kh+t_0$ is described by $g_k=(s_{p_k};q_{k,1},\ldots,q_{k,N+1};R_k)$, and the corresponding configuration manifold is $\G=\Re\times (\Re^3)^{N+1}\times \SO$. This is a Lie group where the group acts on itself by the diagonal action~\cite{MarRat99}: the group action on $s_{p_k}$ and $q_{1,k}\cdots q_{N+1,k}$ is addition, and the group action on $R_k$ is matrix multiplication.

We define a discrete-time kinematics equation using the group action as follows. Define $f_k=(\Delta s_{p_k}; \Delta q_{k,1},\ldots$, $\Delta q_{k,N+1};F_k)\in\G$ such that $g_{k+1}=g_k f_k$:
\begin{align}
    &(s_{p_{k+1}};q_{k+1,1},\ldots,q_{k+1,N+1},R_{k+1})=\nonumber\\
    &(s_{p_k}+\Delta s_{p_k};q_{k,1}+\Delta q_{k,1},\ldots,q_{k,N+1}+\Delta q_{k,N+1};R_k F_k).
\end{align}
Therefore, $f_k\in \G$ represents the relative update between two integration steps. This ensures that the structure of the Lie group configuration manifold is numerically preserved since $g_{k}$ is updated by $f_k$ using the right Lie group action of $\G$ on itself.

\subsection{Discrete Lagrangian}

A discrete Lagrangian $L_d(g_k,f_k):\G\times\G\rightarrow\Re$ is an approximation of the Jacobi solution of the Hamilton--Jacobi equation, which is given by the integral of the Lagrangian along the exact solution of the Euler-Lagrange equations over a single time step:
\begin{align*}
    L_d(g_k,f_k)\approx \int_0^h L(\tilde g(t),{\tilde g}^{-1}(t)\dot{\tilde g} (t))\,dt,
\end{align*}
where $\tilde g(t):[0,h]\rightarrow \G$ satisfies Euler-Lagrange equations with boundary conditions $\tilde{g}(0)=g_k$, $\tilde{g}(h)=g_kf_k$. The resulting discrete-time Lagrangian system, referred to as a variational integrator, approximates the Euler-Lagrange equations to the same order of accuracy as the discrete Lagrangian approximates the Jacobi solution~\cite{Lee08,MarWesAN01}.

We construct a discrete Lagrangian for the string pendulum using the trapezoidal rule. We first find the contributions of each component to the kinetic energy. From the given discretized string pendulum model and \refeqn{Tr}, the kinetic energy of the reel mechanism is approximated by
\begin{align}
T_{k,r} = \frac{1}{2h^2} (\mub s_{p_k} + \kappa_d) \Delta s_{p_k}^2.\label{eqn:Tkr}
\end{align}
Using the chain rule, the partial derivative of $r_k(\sb)$ given by \refeqn{rksb} with respect to $t$ is given by
\begin{align*}
\dot r_k(s) & =
\frac{1}{h} \bigg\{S_0(\zeta_{k,a}) \Delta q_{k,a} + S_1(\zeta_{k,a}) \Delta q_{k,a+1}\\
&\quad+\frac{(L-s)}{(L-s_{p_k})}\frac{(q_{k,a}-q_{k,a+1})}{l_k}\Delta s_{p_k}\bigg\}.
\end{align*}
Substituting this into \refeqn{Ts}, the contribution of the $a$-th string element to the kinetic energy of the string is given by
\begin{align}
T_{k,a} 
& = \frac{1}{2h^2}M^1_k \Delta q_{k,a}\cdot \Delta q_{k,a}
+\frac{1}{2h^2}M^2_k \Delta q_{k,a+1}\cdot \Delta q_{k,a+1}\nonumber\\
&+\frac{1}{2h^2}M^3_{k,a} \Delta s_{p_k}^2
+\frac{1}{h^2}M^{12}_k \Delta q_{k,a}\cdot \Delta q_{k,a+1}\nonumber\\
&+\frac{1}{h^2} M^{23}_{k,a}\Delta s_{p_k} \cdot \Delta q_{k,a+1}
+\frac{1}{h^2} M^{31}_{k,a}\Delta s_{p_k} \cdot \Delta q_{k,a},\label{eqn:Tka}
\end{align}
where inertia matrices are defined in Appendix \ref{sec:IM}.
From the attitude kinetics equations \refeqn{AK}, the angular velocity is approximated by
\begin{align*}
\hat\Omega_k \approx \frac{1}{h} R_k^T (R_{k+1}-R_k) = \frac{1}{h} (F_k -I).
\end{align*}
Define a non-standard inertia matrix $J_d = \frac{1}{2}\mathrm{tr}[J]I-J$. Using the trace operation and the non-standard inertia matrix, the last term of the kinetic energy of the rigid body given by \refeqn{Tr} can be written in terms of $\hat\Omega$ as $\frac{1}{2}\Omega\cdot J\Omega=\frac{1}{2}\mathrm{tr}[\hat\Omega J_d \hat\Omega^T]$. Then, the kinetic energy of the rigid body is given by
\begin{align}
T_{k,b} &= \frac{1}{2h^2} M \Delta q_{k,N+1}\cdot \Delta q_{k,N+1}+\frac{1}{h^2}\tr{(I-F_k)J_d}\nonumber\\
&+\frac{1}{h^2}M \Delta q_{k,N+1}\cdot R_k(F_k-I)\rho_c,\label{eqn:Tkb}
\end{align}
where we use properties of the trace operator: $\mathrm{tr}[AB]=\mathrm{tr}[BA]=\mathrm{tr}[A^TB^T]$ for any matrices $A,B\in\Re^{3\times 3}$. 

From \refeqn{Tkr}, \refeqn{Tka}, \refeqn{Tkb}, the total kinetic energy of the discretized string pendulum is given by
\begin{align}
T_k = T_{k,r} + \sum_{a=1}^N T_{k,a} + T_{k,b}.\label{eqn:Tk}
\end{align}
Similarly, the total potential energy for the given discretized string pendulum can be written as
\begin{align}
V_k & = -\mu g \braces{ (s_{p_k}-d)\, r_d\cdot e_3 
+d^2 (\cos ((s_p-b)/d)-1) }\nonumber\\
&\quad +\sum_{a=1}^N -\frac{1}{2}\mu g l_k e_3 \cdot  (2r_p+q_{k,a}+q_{k,a+1})\nonumber\\
&\quad+\frac{1}{2}\frac{EA}{l_k} (\|q_{k,a+1}-q_{k,a}\| - l_k)^2\nonumber\\
&\quad - M g e_3 \cdot (q_{k,N+1}+r_p+R_k\rho_c). \label{eqn:Vk}
\end{align}

This yields the discrete-Lagrangian of the discretized string pendulum
\begin{align}
L_{d_k}(g_k,f_k) = h T_k(g_k,f_k) -\frac{h}{2} V_k(g_k,f_k) - \frac{h}{2} V_{k+1}(g_k,f_k).\label{eqn:Ld}
\end{align}

\subsection{Lie Group Variational Integrator}

For a given discrete Lagrangian, the discrete action sum is given by $\mathfrak{G}_d=\sum_{k} L_{d_k}$. As the discrete Lagrangian approximates the action integral over a single discrete time step, the action sum approximates the action integral. According to the discrete Lagrange--d'Alembert principle, the sum of the variation of the action sum and the discrete virtual work done by external control moments and constraints is equal to zero. This yields discrete-time forced Euler-Lagrange equations referred to as variational integrators~\cite{MarWesAN01}. This procedure is followed for an arbitrary discrete Lagrangian defined on a Lie group configuration manifold in~\cite{Lee08} to obtain a Lie group variational integrator:
\begin{gather}
\begin{aligned}
    \T_e^*\L_{f_{k-1}}\cdot  &\D_{f_{k-1}} L_{d_{k-1}}-\Ad^*_{f_{k}^{-1}}\cdot(\T_e^*\L_{f_{k}}\cdot \D_{f_{k}}L_{d_{k}})\\
    &+\T_e^*\L_{g_{k}}\cdot \D_{g_{k}} L_{d_{k}}+u_{d_k}+Q_{d_k}=0,
\end{aligned}\label{eqn:DEL_G0}\\
    g_{k+1} = g_k f_k,\label{eqn:DEL_G1}
\end{gather}
where $\T^*\L:\G\times\T^*\G\rightarrow\T^*\G$ is the co-tangent lift of the left translation action, $\D_f$ represents the derivative with respect to $f$, and $\Ad^*:\G\times\g^*\rightarrow\g^*$ is the co-Adjoint operator~\cite{MarRat99}. 

The contribution of the external control moment and the Carnot energy loss term are denoted by $u_{d_k}$ and $Q_{d_k}$. They are defined to approximate the additional term in the variational principle \refeqn{VP} that arises due to a discontinuity:
\begin{align*}
\int_{kh}^{(k+1)h} (Q+u/d)\delta s_p \,dt \approx (Q_{d_k} + u_{d_k}) \delta s_{p_k}.
\end{align*}
From \refeqn{dotE}, these are chosen as
\begin{align}
Q_{d_k} &= -\frac{h}{2l_k^2} (\mub \Delta s_{p_k}^2 /h^2 +EA) ( \norm{q_{k,2}}-l_k)^2,\label{eqn:Qdk}\\
u_{d_k} &= h u_k /d.\label{eqn:udk}
\end{align}

We substitute the expressions for the discrete Lagrangian \refeqn{Ld}, the Carnot energy loss term \refeqn{Qdk}, and the control moment \refeqn{udk} into \refeqn{DEL_G0} and \refeqn{DEL_G1} to obtain a Lie group variational integrator for the discretized string pendulum model. This involves deriving the derivatives of the discrete Lagrangian and their co-tangent lift. The detailed procedure and the resulting expressions for Lie group variational integrators are summarized in the Appendix.

\paragraph{Computational Properties}

The proposed Lie group variational integrators have desirable computational properties. The Lie group configuration manifold is often parameterized. But, the local parameterizations of the special orthogonal group, such as Euler angles or Rodrigues parameters, have singularities. In a numerical simulation of large angle maneuvers of a string pendulum, these local parameters should be successively switched from one type to another in order to avoid their singularities. They also lead to excessive complexity.
 Non-parametric representations such as quaternions also have associated difficulties: there is an ambiguity in representing an attitude since the group $\mathsf{SU}(2)$ of quaternions double cover $\SO$. Furthermore, as the unit length of quaternions is not preserved in numerical simulations, attitudes cannot be determined accurately. Sometimes, numerical solutions updated by any one-step integration method are projected onto the Lie group at each time step~\cite{DieRusSJNA94}. Such projections may destroy the desirable long-time behavior of one-step methods, since the projection typically corrupts the numerical results.  Lie group variational integrators update the group elements by using a group operation. Therefore, the Lie group structure is naturally preserved at the level of machine precision, and they avoid any singularity and complexity associated with other approaches.

Since Lie group variational integrators are constructed according to Hamilton's principle, their numerical trajectories preserve a symplectic form and a momentum map associated with any symmetry. These ensure long-term structural stability and avoid 
artificial numerical dissipation. These properties are difficult to achieve in conventional approaches based on finite difference approximation of continuous equations of motion. 

In summary, the proposed Lie group variational integrators for a string pendulum will be particularly useful when studying nontrivial maneuvers that combine large elastic deformations and large rigid motions accurately over a long time period.

\section{Numerical Examples}\label{sec:NE}

We now numerically demonstrate the computational properties of the Lie group variational integrators developed in the previous section. The properties of the reel mechanism are as follows: $b=d=0.5\,\mathrm{m}$, $\kappa_d=1\,\mathrm{kg}$. The material properties of the string are chosen to represent a rubber string~\cite{KuhSeiAMC95}: $\mub=0.025\,\mathrm{kg/m}$, $EA=40\,\mathrm{N}$, $L=100\,\mathrm{m}$. The rigid body is chosen as an elliptic cylinder with a semimajor axis $0.5\,\mathrm{m}$, a semiminor axis $0.4\,\mathrm{m}$, and a height $0.8\,\mathrm{m}$. The mass and the location of the center of mass of the rigid body are $M=0.1\,\mathrm{kg}$, and $\rho_c=[0.3,0.2,0.4]\,\mathrm{m}$, respectively.

We consider three cases: (1) dynamics of a fixed length string pendulum released from a horizontal configuration, (2) deployment dynamics due to gravity from a horizontal configuration, and (3) retrieval dynamics using a constant control moment. Initial conditions are as follows:

{\small\selectfont
\begin{center}
\begin{tabular}{cc|ccc}
& $s_{p_0}\,\mathrm{(m)}$ & $q_{0,a}\,\mathrm{(m)}$ & $u_k\,\mathrm{(Nm)}$\\\hline
(1)  & 90 & $l_0 e_1$ & -\\
(2)  & 99 & $l_0(a-1)e_1$ & 0\\
(3)  & 90 & $l_0(a-1)\,(\sin 15^\circ e_1+\cos 15^\circ e_3)$ & 2.09\\
\end{tabular}
\end{center}}

For all cases, we choose $\dot s_p = 0\,\mathrm{m/s}$, $\dot q_{0,a}=0\,\mathrm{m/s}$ for $a\in\{1,N\}$, $\dot q_{0,N+1}=0.5e_2\,\mathrm{m/s}$,  $R_0=I$, $\Omega_0=0\,\mathrm{rad/sec}$. The deployed portion of the string is discretized by $N=20$ elements, and the time step is $h=0.0005$ second. Simulation time is $T=10$, $T=8$, and $T=10$ seconds for each case, respectively.

\paragraph{Energy transfer} 

The following figures show the simulation results. The maneuver of the string pendulum is illustrated by snapshots, where the relative elastic potential distribution at each instant is denoted by color shading  (the corresponding animations are available at http://my.fit.edu /\~{}taeyoung). As the point where the string is attached to the rigid body is displaced from the center of mass of the rigid body, the rigid body dynamics are directly coupled to the elastic string dynamics. The illustrated maneuvers clearly show the nontrivial coupling between the strain deformation, the rigid body dynamics, and the reel mechanism. 

\begin{figure*}
\vspace*{3.5cm}

\centerline{
	\hspace*{0.05\textwidth}		
	\subfigure[Snapshots at each 0.2 second $t\in [ 0,5{]} $]{
		\includegraphics[width=0.40\textwidth]{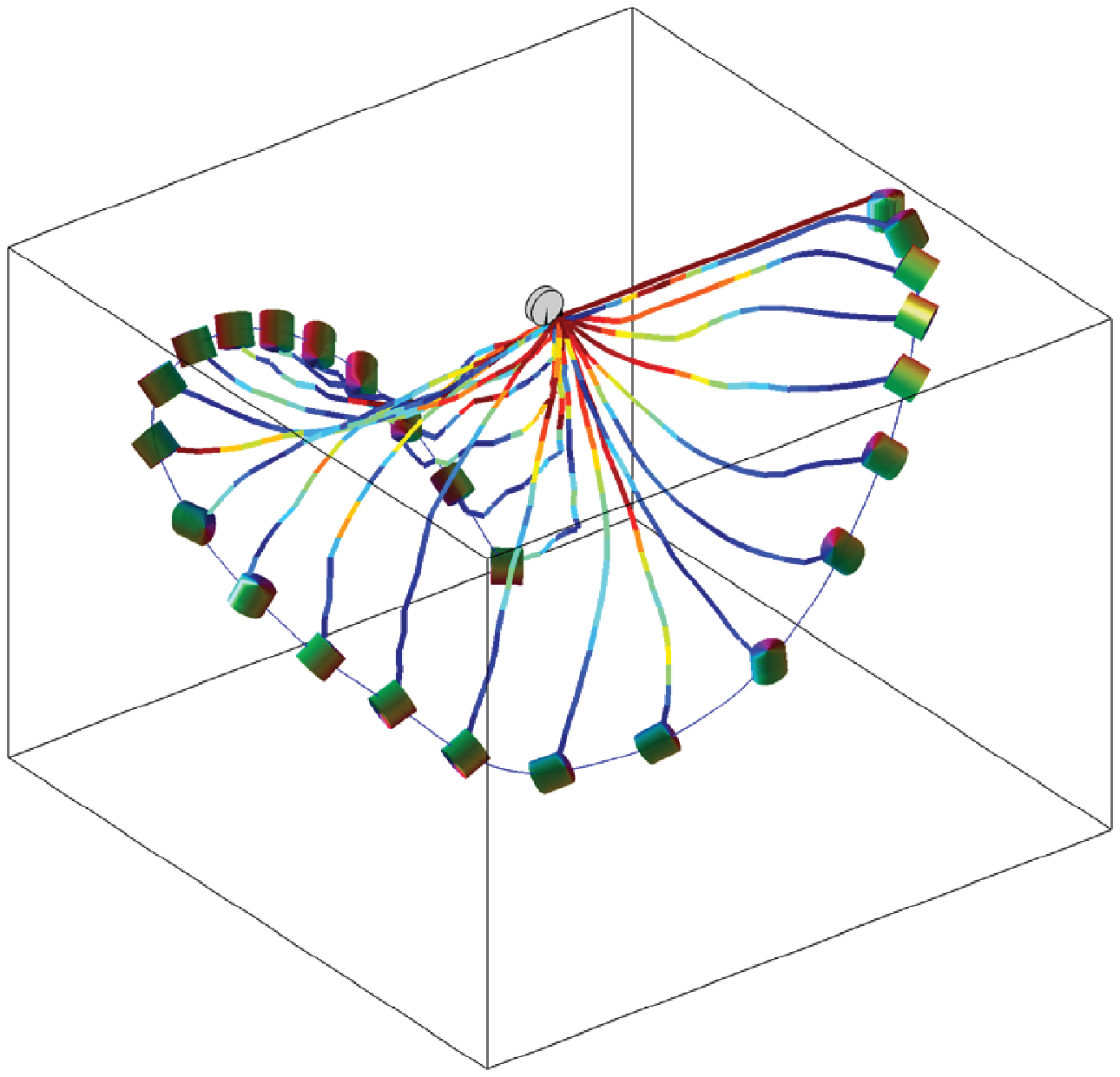}\label{fig:SS1}}
	\hspace*{0.08\textwidth}		
	\subfigure[Energy exchange ($E$:solid, $T$:solid, $T_{rot}$:dashed, $V_{gravity}$:dash-dotted, $V_{elastic}$:dotted)]{
\renewcommand{\xyWARMinclude}[1]{\includegraphics[width=0.48\textwidth]{#1}}%616
{\footnotesize\selectfont
$$\begin{xy}
\xyWARMprocessEPS{Lee2b}{eps}%%
\xyMarkedImport{}%%
\xyMarkedMathPoints{1-15}
\end{xy}
$$}\label{fig:ET1}}
}
\centerline{
	\subfigure[Stretched length of the deployed portion of the string, and the second component of the angular velocity $\Omega$]{
		\includegraphics[width=0.48\textwidth]{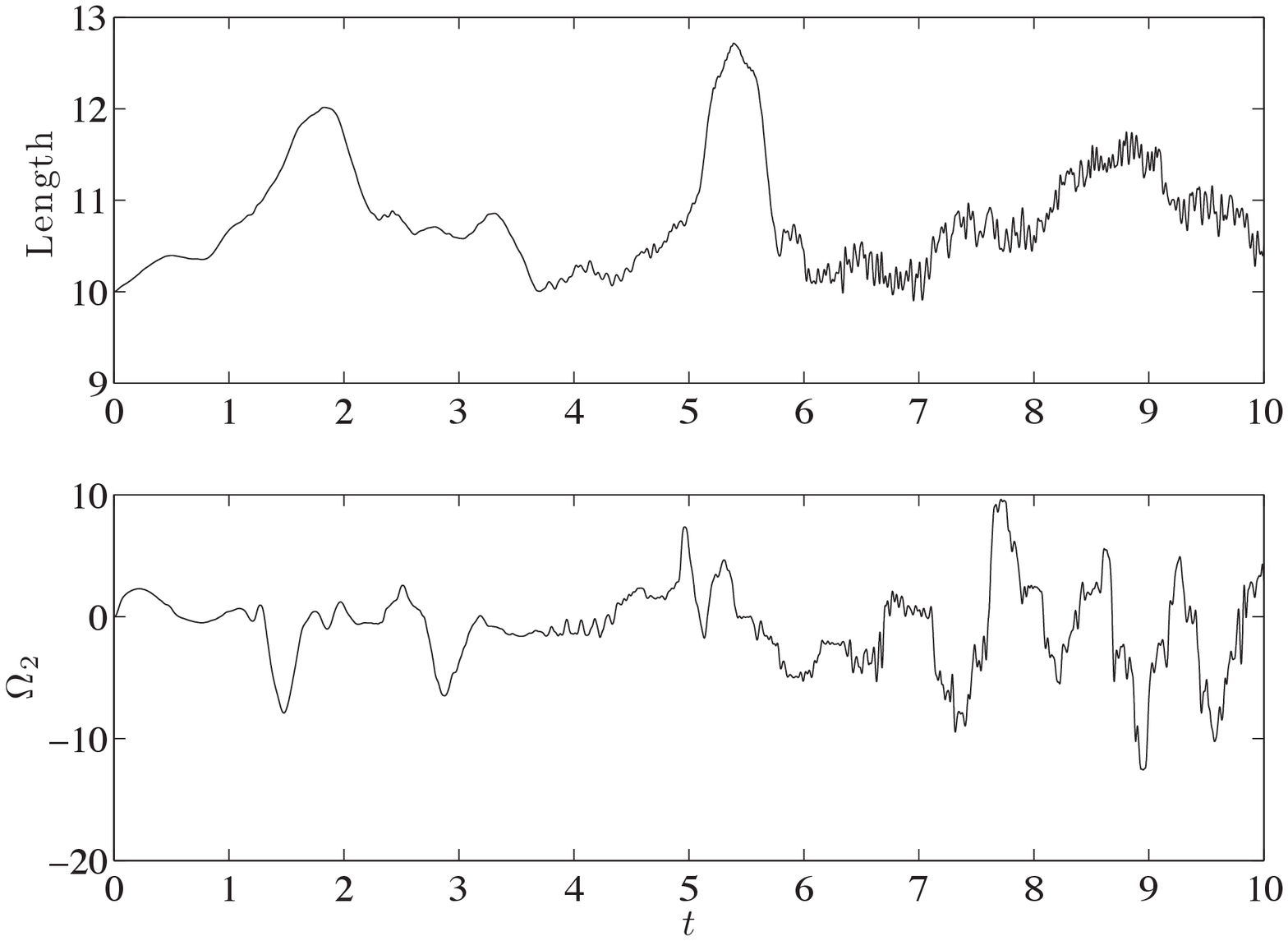}}
	\hspace*{0.05\textwidth}		
	\subfigure[Deviation of conserved quantities: total energy and the total angular momentum about the gravity direction]{
		\includegraphics[width=0.47\textwidth]{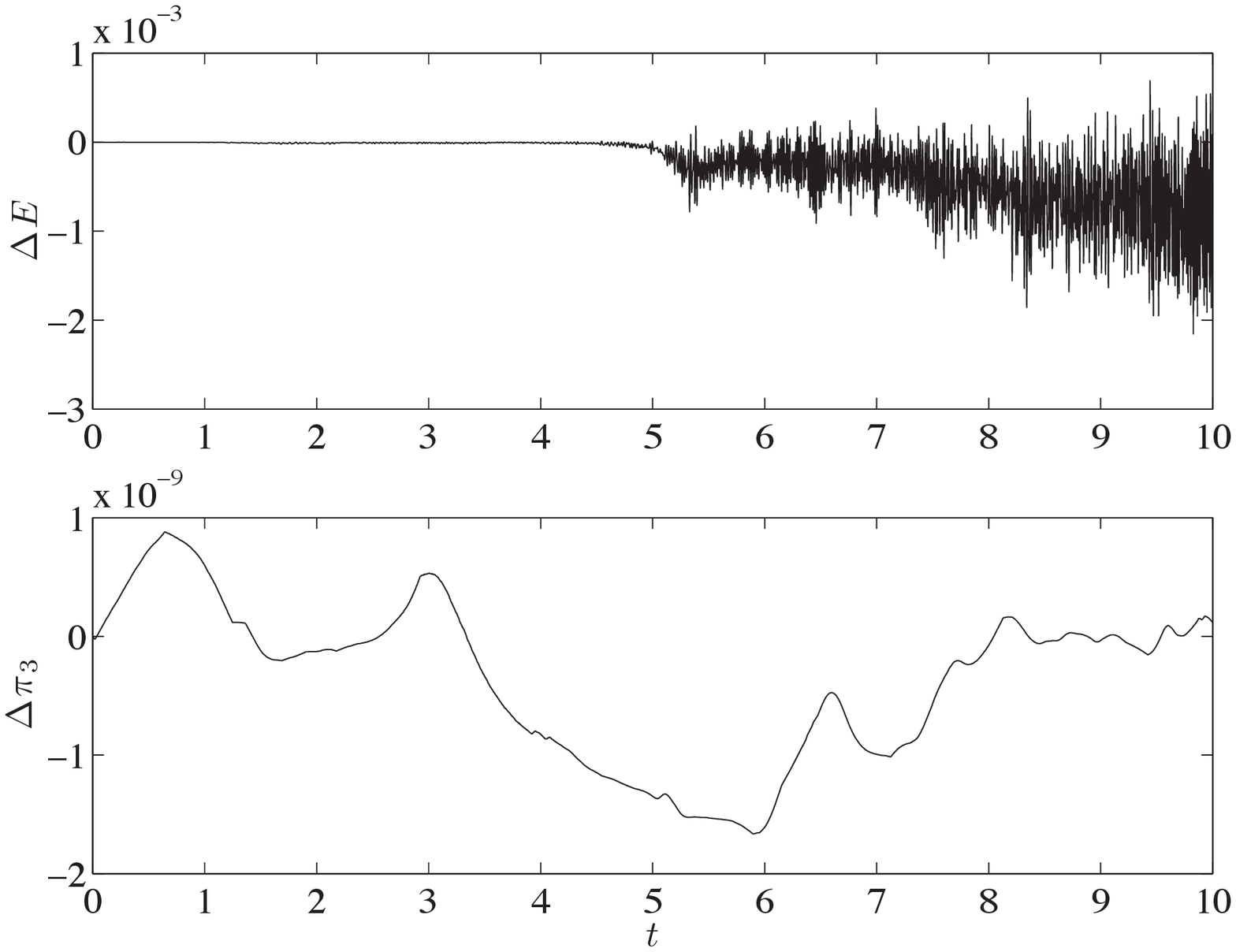}\label{fig:1conv}}
}
\caption{Fixed length string pendulum}
\end{figure*}

\begin{figure*}
\vspace*{3.5cm}

\centerline{
	\hspace*{0.07\textwidth}		
	\subfigure[Snapshots at each 0.4 second $t\in [ 0,8{]} $]{
		\includegraphics[width=0.34\textwidth]{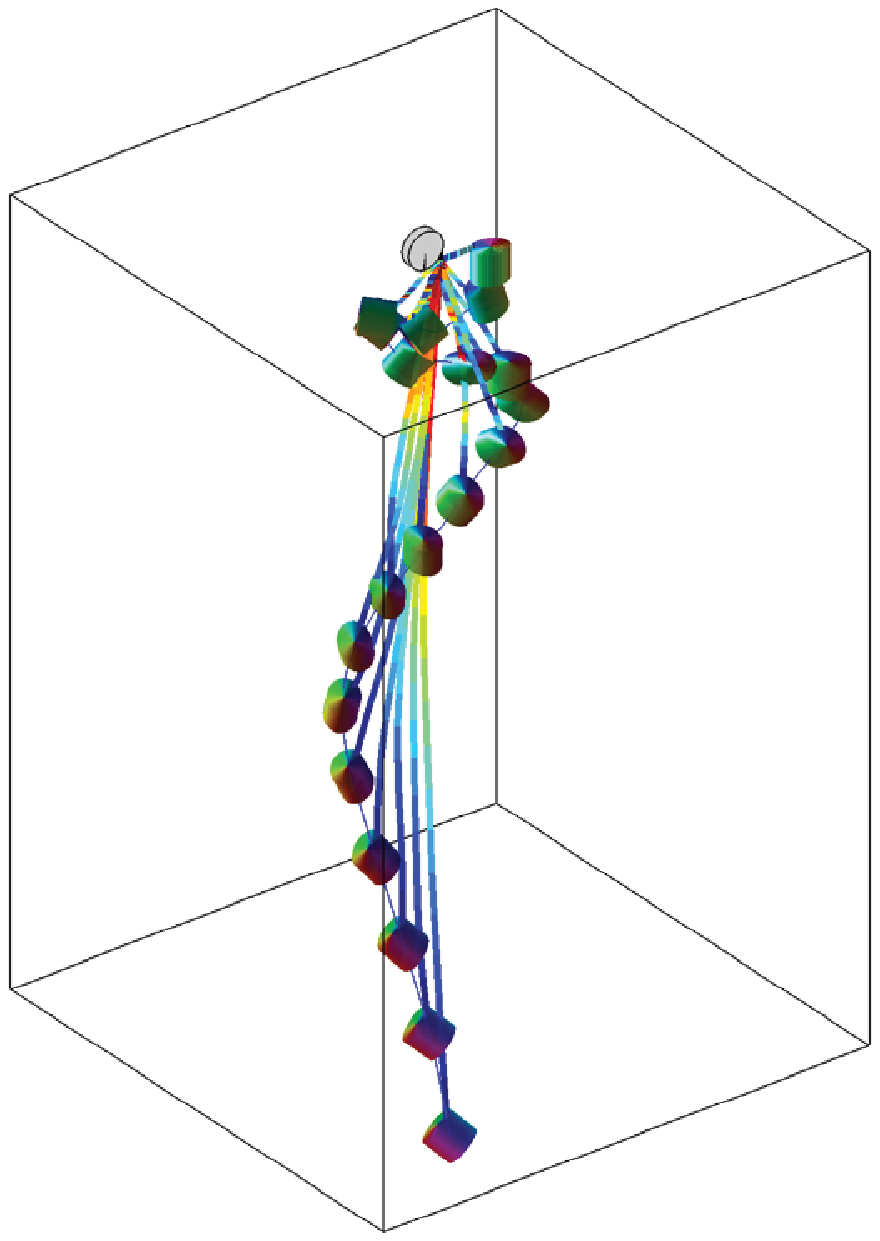}}
	\hspace*{0.10\textwidth}		
	\subfigure[Energy exchange ($E$:solid, $T$:solid, $T_{rot}$:dashed, $V_{gravity}$:dash-dotted, $V_{elastic}$:dotted)]{
\renewcommand{\xyWARMinclude}[1]{\includegraphics[width=0.47\textwidth]{#1}}
{\footnotesize\selectfont
$$\begin{xy}
\xyWARMprocessEPS{Lee3b}{eps}
\xyMarkedImport{}
\xyMarkedMathPoints{1-15}
\end{xy}
$$}\label{fig:ET2}}
}
\centerline{
	\subfigure[Length of the deployed portion of the string (stretched:solid, unstreched:dashed), and the second component of the angular velocity $\Omega$]{
		\includegraphics[width=0.47\textwidth]{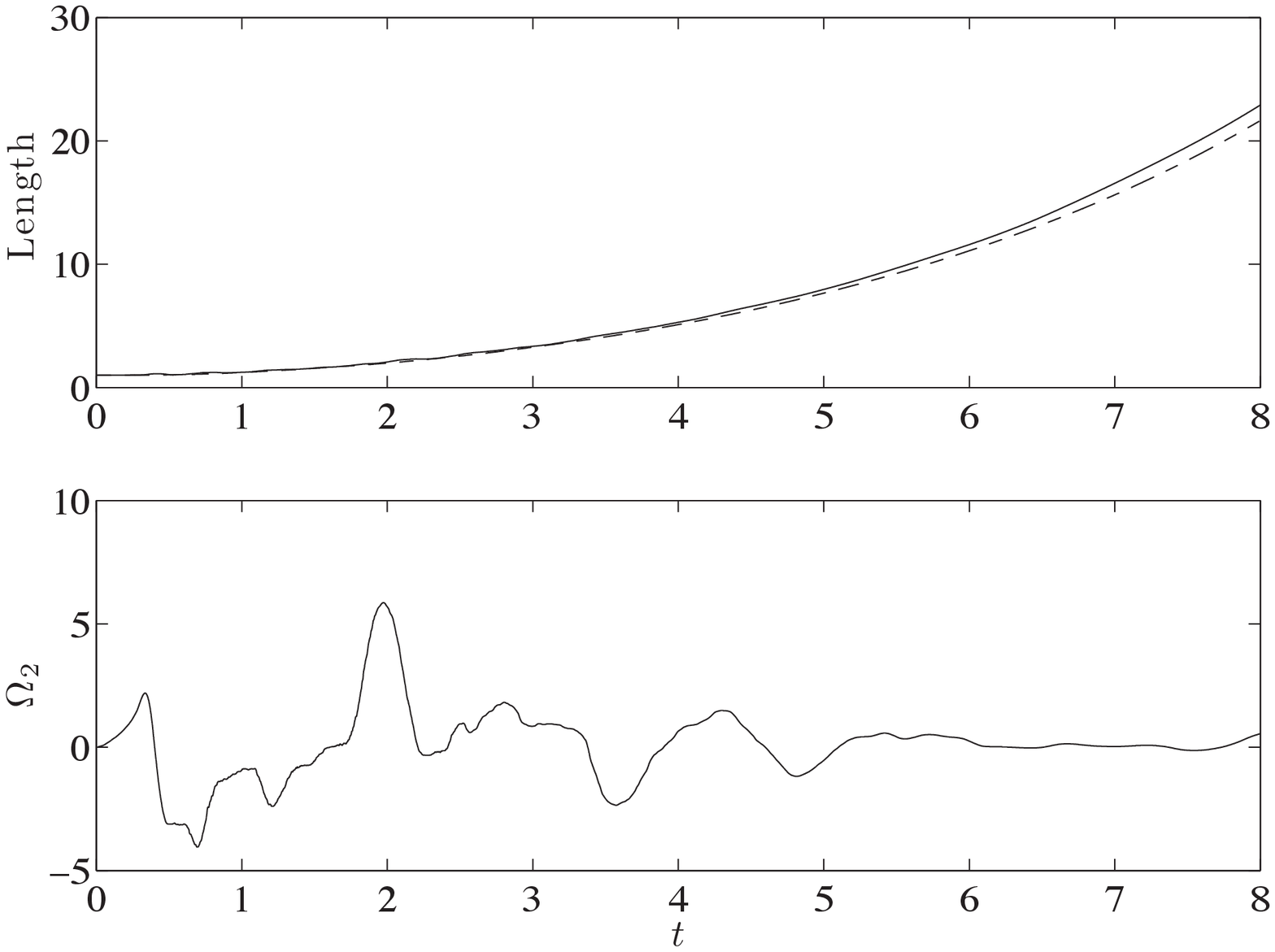}}
	\hspace*{0.04\textwidth}		
	\subfigure[Deviation of conserved quantities: the difference between the computed total energy change and the energy dissipation due to the velocity discontinuity, the orthogonality error of the rotation matrix]{
		\includegraphics[width=0.48\textwidth]{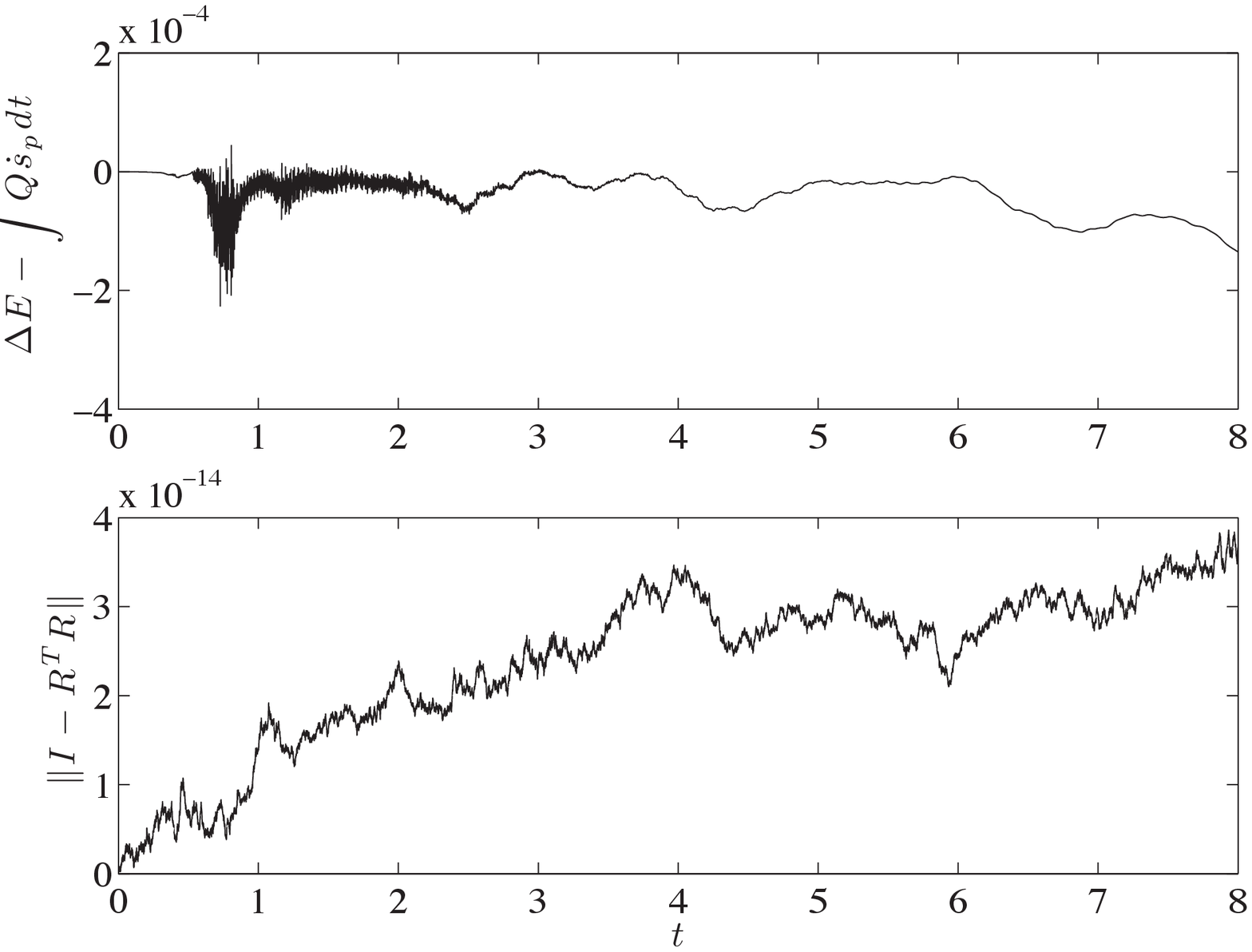}\label{fig:2conv}}
}
\caption{Deployment due to gravity}
\end{figure*}

The energy exchange plots also show that there is significant energy transfer between the kinetic energy, the gravitational potential energy, and the elastic potential energy. In \reffig{ET1}, there is an energy exchange between the kinetic energy and the gravitational potential energy. But, when the string is mostly stretched at $t=1.8$ and $t=5.5$ seconds, part of the kinetic energy is transferred to the elastic potential energy and the rotational kinetic energy. The rigid body starts to tumble at $t=7$ seconds. The elastic potential energy transfer along the string is observed in \reffig{SS1}. For the deployment case shown in \reffig{ET2}, the gravitational potential energy is generally transferred to the kinetic energy. As the length of the deployed portion of the string increases, the elastic potential energy increases and the rotational kinetic energy decreases. For the third retrieval case, both the gravitational potential energy and the total energy increase due to the constant control moment. In addition, there is a smaller-scale periodic energy exchange between the elastic potential and the gravitational potential energy with an approximate period of 1.2 seconds. The rigid body starts tumbling at $t=3$ seconds.

\paragraph{Conservation Properties}

\begin{figure*}
\vspace*{3.5cm}

\centerline{
	\hspace*{0.02\textwidth}		
	\subfigure[Snapshots at each 0.5 second $t\in [ 0,10{]} $]{
		\includegraphics[width=0.43\textwidth]{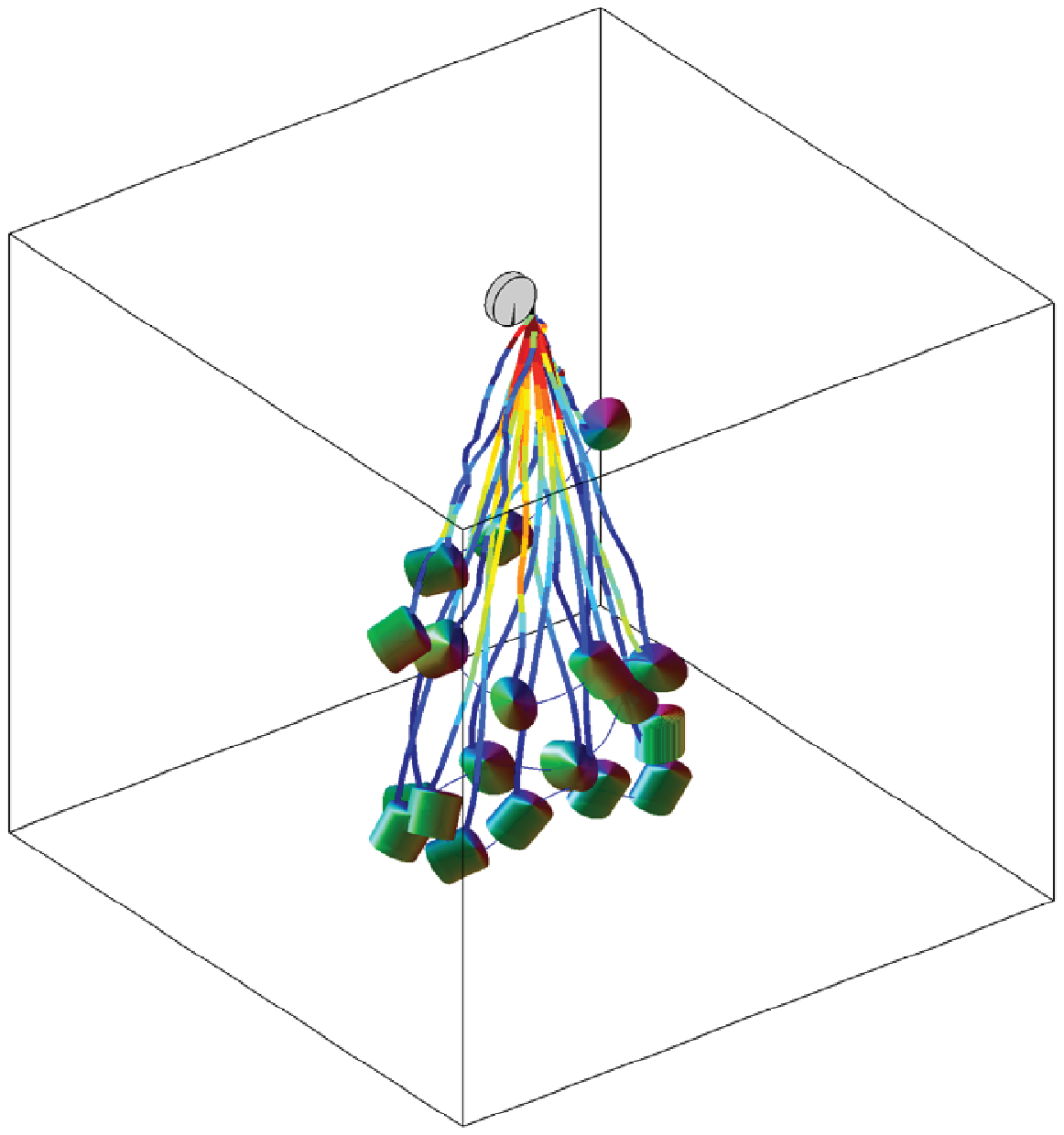}}
	\hspace*{0.05\textwidth}		
	\subfigure[Energy exchange ($E$:solid, $T$:solid, $T_{rot}$:dashed, $V_{gravity}$:dash-dotted, $V_{elastic}$:dotted)]{
\renewcommand{\xyWARMinclude}[1]{\includegraphics[width=0.47\textwidth]{#1}}
{\footnotesize\selectfont
$$\begin{xy}
\xyWARMprocessEPS{Lee4b}{eps}
\xyMarkedImport{}
\xyMarkedMathPoints{1-15}
\end{xy}
$$}\label{fig:ET3}}
}
\centerline{
	\subfigure[Length of the deployed portion of the string (stretched:solid, unstreched:dashed), and the second component of the angular velocity $\Omega$]{
		\includegraphics[width=0.47\textwidth]{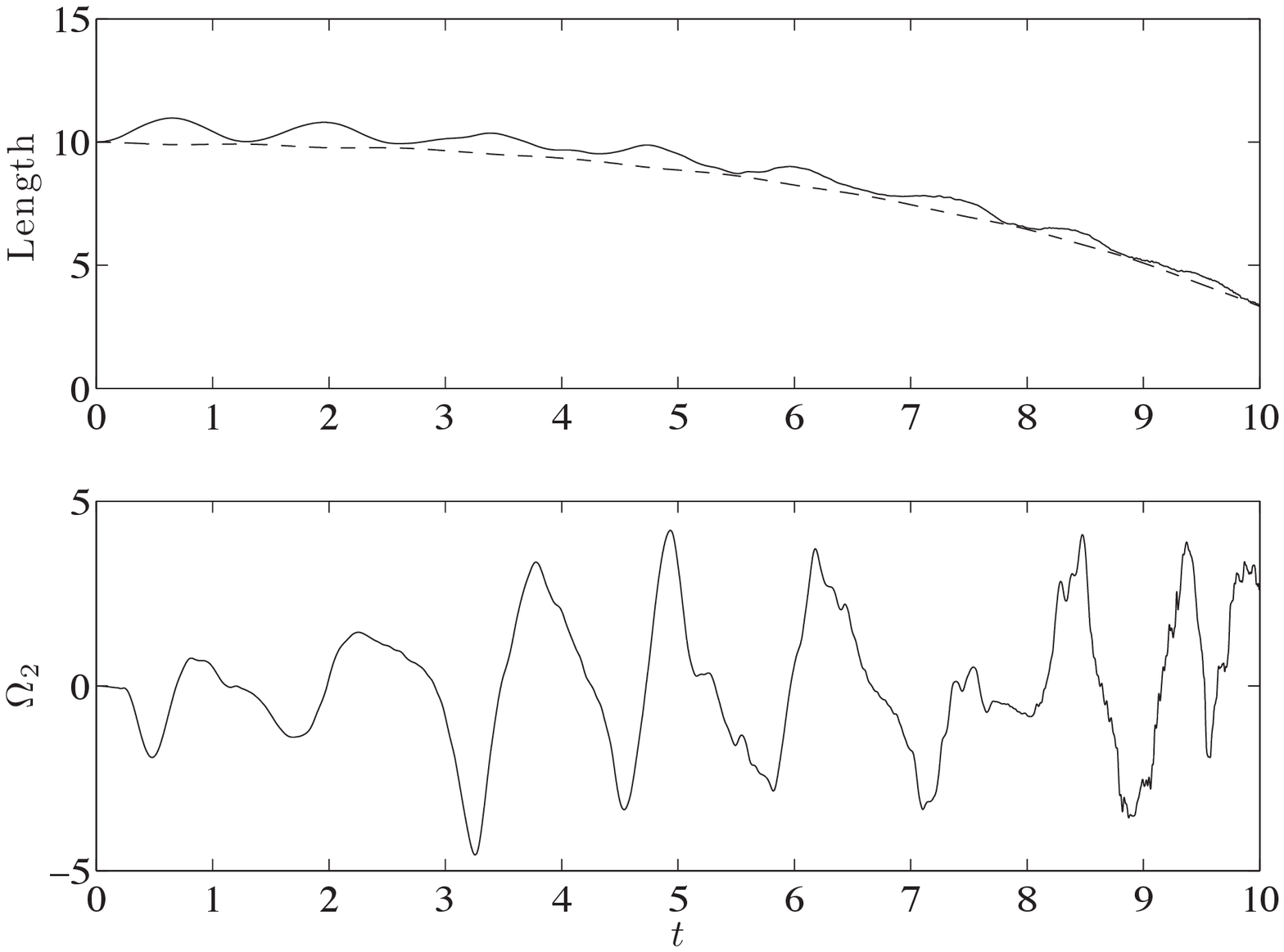}}
	\hspace*{0.04\textwidth}		
	\subfigure[Orthogonality error of the rotation matrix]{
		\includegraphics[width=0.48\textwidth]{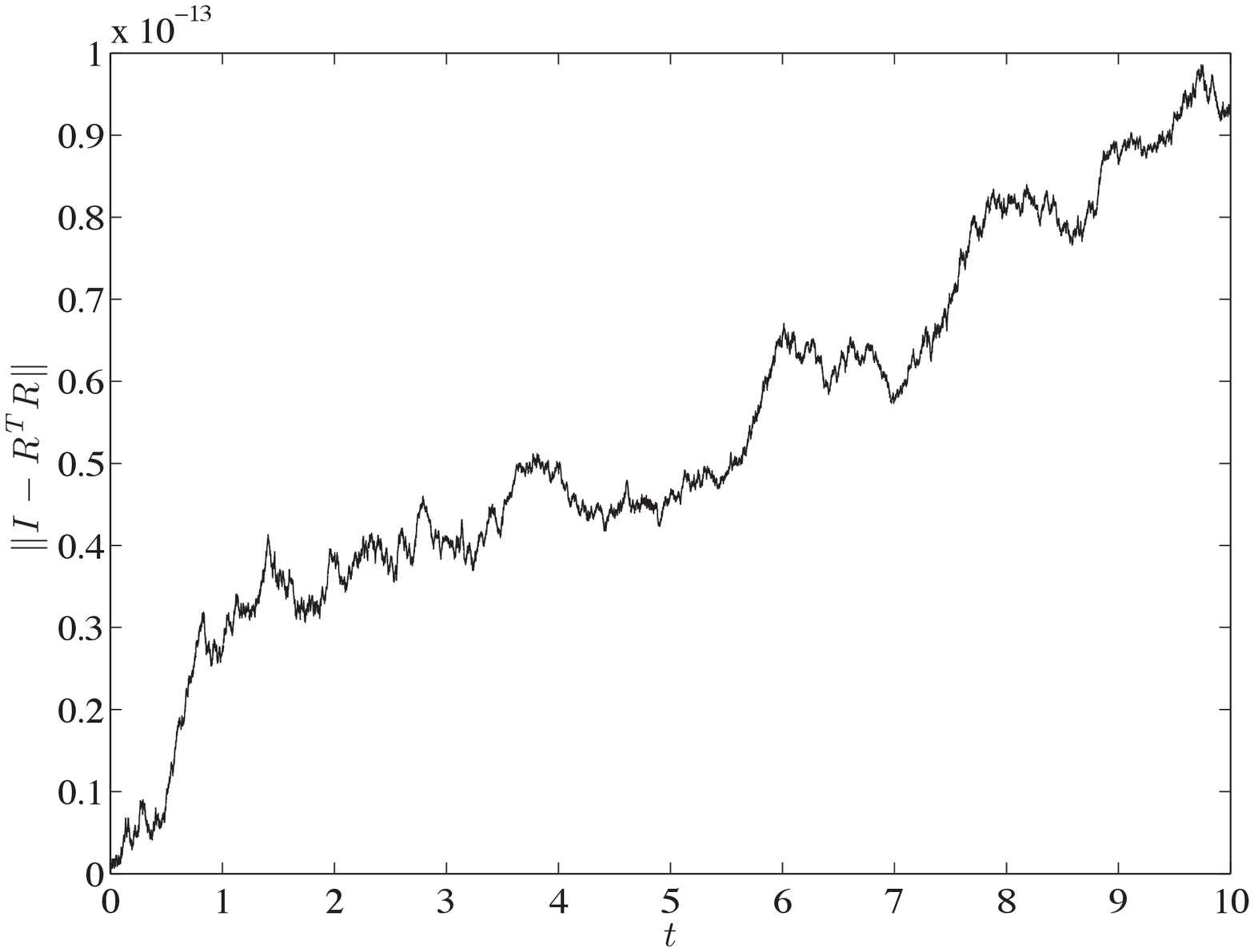}\label{fig:3conv}}
}
\caption{Retrieval using a constant control moment}

\end{figure*}

The proposed Lie group variational integrators exhibit excellent conservation properties for these complicated maneuvers of the string pendulum. For the fixed length string dynamics, the total energy and the total angular momentum about the gravity direction should be preserved. The deviations of those quantities are shown in \reffig{1conv}, where the maximum deviation of the total energy is less than $0.01\%$ of the maximum kinetic energy, and the deviation of the angular momentum is less than 
$3\times 10^{-8}\%$ of its initial value. For the second deployment case, the total energy dissipates only due to the velocity discontinuity. \reffig{2conv} shows the difference between the computed total energy change and the energy dissipation computed by the Carnot energy loss term \refeqn{dotE}. The difference is less than $0.0003\%$ of the maximum kinetic energy, which illustrates that there is no artificial numerical dissipation caused by the proposed Lie group variational integrator. 
The orthogonal structure of rotation matrices is preserved to machine precision.   \reffig{2conv} and \ref{fig:3conv} show that the orthogonality error, measured by $\norm{I-R^T R}$, is less than $10^{-13}$.

\section{Conclusions}\label{sec:CON}

We have developed continuous-time equations of motion and geometric numerical integrators, referred to as Lie group variational integrators, for a 3D elastic string pendulum attached to a rigid body and a reel mechanism. They are carefully derived while taking account of the length change of the deployed portion of the string, the Lie group configuration manifold of the rigid body, and the velocity discontinuity at the guide way entrance. The continuous-time equations of motion provide an analytical model that is defined globally on the Lie group configuration manifold.  The Lie group variational integrator preserves the geometric features of the system, thereby yielding a reliable numerical method to compute the nonlinear coupling between the large string deformation and the nontrivial rigid body dynamics accurately over a long time period. In short, this paper provides high fidelity analytical and computational models for a string pendulum.

The numerical experiments suggest that accurately modeling the reeling mechanism is of critical importance in order to capture the correct dynamics, due to the disturbance that is introduced in the string at the point of contact with the reeling mechanism when the string is deployed or retracted. One can observe that this disturbance propagates down the string at a velocity that is determined by the elastic properties of the string. Since the point of contact between the string and the rigid body does not go through the center of mass of the rigid body, the elastic disturbance excites a rotational response in the rigid body. As such, accurately modeling the reel mechanism, elastic string dynamics, rigid body motion, and their interactions, is critical for obtaining realistic predictions about how towed underwater vehicles and tethered spacecraft behave when performing aggressive maneuvers.

The proposed string pendulum model and computational approach can be extended in several ways. For example, different types of string models can be considered, such as an inextensible string, nonlinear elasticity, and bending stiffness. The reel mechanism can be generalized by assuming that the portion of the string on the drum is also extensible.
These results can be extended to model tethered spacecraft in orbit, and they can be used to study associated optimal control problems by adopting the \textit{discrete mechanics and optimal control} approach~\cite{JuMaOb2005}.

\clearpage\newpage

\clearpage\newpage

\clearpage\newpage

\clearpage\newpage 

\appendix

\section{Development of the Lie Group Variational Integrator for a String Pendulum}

\subsection{Inertia Matrices for the Discrete Lagrangian}\label{sec:IM}
The inertia matrices for the discrete Lagrangian are defined as follows.
\begin{align*}
M^1_k&=\frac{1}{3}\mub l_k,\qquad 
M^2_k=M^1_k,\\
M^3_{k,a} &= \frac{1}{3}\mub l_k \frac{(3N^2+3N+1-6Na-3a+3a^2)}{N^2},
\\
M^{12}_k &= \frac{1}{6}\mub l_k,\qquad
M^{23}_{k,a} = \frac{1}{6}\mub \frac{(1+3N-3a)}{N}(q_{k,a}-q_{k,a+1}),\\
M^{31}_{k,a} &= \frac{1}{6}\mub \frac{(2+3N-3a)}{N}(q_{k,a}-q_{k,a+1}).
\end{align*}

\subsection{Derivatives of the discrete Lagrangian}\label{sec:DDL}

The Lie group variational integrator given by \refeqn{DEL_G0} is expressed in terms of the derivatives of the discrete Lagrangian and their co-tangent lift. Here, we describe how to compute the co-tangent lift and the co-Adjoint operator on the configuration manifold $\G=\Re\times(\Re^3)^{N+1}\times\SO$ without introducing the formal definition of those operators.

The co-tangent lift of the left translation on a real space is the identity map on that real space. Using the product structure of the configuration manifold $\G=\Re\times(\Re^3)^{N+1}\times\SO$, the derivative of the discrete Lagrangian with respect to $f_k=(\Delta s_{p_k}; \Delta q_{k,1},\ldots,\Delta q_{k,N+1};F_k)\in\G$ is given by
\begin{align}
\T_e^*\L_{f_{k}}\cdot \D_{f_{k}}L_{d_{k}} = \big[& \D_{\Delta s_{p_k}} L_{d_k};\, \D_{\Delta q_{1,k}} L_{d_k},\cdots,\nonumber\\
 &\D_{\Delta q_{N+1,k}} L_{d_k};\, \T_I^*\L_{F_{k}}\cdot \D_{F_{k}}L_{d_{k}}\big].
\end{align}

Deriving the derivatives of the discrete Lagrangian with respect to $\Delta s_{p_k}$ or $\Delta q_{k,a}$ is straightforward. For example, from \refeqn{Tka}, \refeqn{Tk}, \refeqn{Ld}, the derivative of the discrete Lagrangian with respect to $\Delta q_{k,a}$ for any $a\in\{2,\ldots,N\}$ is given by
\begin{align}
\D_{\Delta q_{k,a}}  L_{d_k} & = h \D_{\Delta q_{k,a}} T_{k,a-1} + h \D_{\Delta q_{k,a}} T_{k,a-1} -\frac{h}{2}\D_{\Delta q_{k,a}} V_{k+1}\nonumber\\
& = \frac{1}{h} M^{12}_k \Delta q_{k,a-1} + \frac{2}{h} M^1_k \Delta q_{k,a} +\frac{1}{h} M^{12}_k \Delta q_{k,a+1}\nonumber\\
& +\frac{1}{h}(M^{31}_{k,a}+M^{23}_{k,a-1})\Delta s_{p_k}
-\frac{h}{2}\D_{q_{k+1,a}} V_{k+1},\label{eqn:DdelqLd}
\end{align}
where the derivative of the potential energy is given by
\begin{align}
\D_{q_{k,a}} V_k & = -\mu g l_k e_3+\nabla V^e_{k,a-1} - \nabla V^e_{k,a},\\
\nabla V^e_{k,a} & =\frac{EA}{l_k} \frac{\norm{q_{k,a+1}-q_{k,a}}-l_k}{\norm{q_{k,a+1}-q_{k,a}}} (q_{k,a+1}-q_{k,a}).
\end{align}
Expressions for the other derivatives of the discrete Lagrangian with respect to $s_{p_k},\Delta s_{p_k},q_{k,a}$ are similarly developed and they are summarized later.

Now we find the derivative of the discrete Lagrangian with respect to $F_k$. From \refeqn{Tkb}, \refeqn{Vk}, \refeqn{Ld}, we have
\begin{align*}
\D_{F_{k}}L_{d_{k}}\cdot \delta F_k
& = \frac{1}{h} \tr{-\delta F_k J_d} + \frac{M}{h} R_k^T \Delta q_{k,N+1}\cdot \delta F_k \rho_c\\
& \quad+\frac{h}{2} Mg e_3\cdot R_k \delta F_k \rho_c.
\end{align*}
Similar to \refeqn{delR}, the variation of the rotation matrix $F_k$ can be written as $\delta F_k = F_k \hat\zeta$ for a vector $\zeta\in\Re^3$. From the definition of the co-tangent lift of the left translation, we have
\begin{align*}
(\T_I^*\L_{F_{k}}\cdot \D_{F_{k}}L_{d_{k}}) \cdot \zeta & = 
\frac{1}{h} \tr{- F_k\hat\zeta_k J_d} + \frac{M}{h} R_k^T \Delta q_{k,N+1}\cdot F_k\hat\zeta_k \rho_c\\
&\quad +\frac{h}{2} Mg e_3\cdot R_k  F_k\hat\zeta_k \rho_c.
\end{align*}
By repeatedly applying the following property of the trace operator, $\mbox{tr}[AB]=\mbox{tr}[BA]=\mbox{tr}[A^TB^T]$ for any $A,B\in\Re^{3\times 3}$, the first term can be written as $\mbox{tr}[-F_{k} \hat\zeta_{k} J_{d}]= \mbox{tr}[- \hat\zeta_{k} J_{d}F_{k}]=\mbox{tr}[\hat\zeta_{k}F_{k}^TJ_{d}]=-\frac{1}{2}\mbox{tr}[\hat\zeta_{k}(J_{d}F_{k}-F_{k}^T J_{d})]$. Using the property of the hat map, $x^T y = -\frac{1}{2}\mbox{tr}[\hat x\hat y]$ for any $x,y\in\Re^3$, this can be further written as $((J_{d}F_{k}-F_{k}^T J_{d})^\vee) \cdot \zeta_{k}$. As $y\cdot\hat x z = \hat z y\cdot x$ for any $x,y,z\in\Re^3$, the second term can be written as $F_k^TR_k^T \Delta q_{k,N+1} \cdot \hat\zeta_k \rho_c=\hat\rho_c F_k^TR_k^T \Delta q_{k,N+1}\cdot \zeta_k$. Using these, we obtain
\begin{align}
\T_I^*\L_{F_{k}}\cdot \D_{F_{k}}L_{d_{k}}& =\frac{1}{h} (J_dF_k -F_k^TJ_d )^\vee + \frac{M}{h} \hat\rho_c F_k^T R_k^T \Delta q_{k,N+1}\nonumber \\&\quad + \frac{h}{2} Mg\hat \rho_c F_k^T R_k^T e_3.\label{eqn:TLDFLd}
\end{align}
Expression for the derivatives of the discrete Lagrangian with respect to $R_k$ is similarly developed.

In summary, in addition to \refeqn{DdelqLd} and \refeqn{TLDFLd}, derivatives of the discrete Lagrangian are summarized as follows.
{\allowdisplaybreaks
\begin{align*}
\D_{\Delta s_{p_k}} L_{d_k} 
& = \frac{1}{h}M^0_k\Delta s_{p_k}
+\frac{1}{h}\sum_{a=2}^{N} (M^{31}_{k,a}+M^{23}_{k,a-1})\cdot \Delta q_{k,a}
\\& + \frac{1}{h} M^{23}_{k,N} \Delta q_{k,N+1} -\frac{h}{2}\D_{s_{p_{k+1}}} V_{k+1},\\
M^0_k & = \mu s_{p_k} + \kappa_d + \frac{1}{3}\mu(L-s_{p_k}),\\
\D_{\Delta q_{k,N+1}} L_{d_k} & = \frac{1}{h} (M^2_k+M) \Delta q_{k,N+1} +\frac{1}{h} M_k^{12} \Delta q_{k,N} \\&+ \frac{1}{h} M^{23}_{k,N} \Delta s_{p_k}
+\frac{1}{h} M R_k(F_k-I)\rho_c 
\\&-\frac{h}{2} \D_{q_{k+1,N+1}} V_{k+1},
\\
\D_{s_{p_k}} L_{d_k} 
& = \frac{\mu}{3 h}  \Delta s_{p_k}^2 
- \frac{\mu}{6N h} \sum_{a=1}^N (\Delta q_{k,a}\cdot \Delta q_{k,a}\\&
+\Delta q_{k,a+1}\cdot \Delta q_{k,a+1}
+\Delta q_{k,a}\cdot \Delta q_{k,a+1})
\\&-\frac{h}{2} \D_{s_{p_k}} V_k - \frac{h}{2} \D_{s_{p_{k+1}}} V_{k+1}
,\\
\D_{q_{k,a}} L_{d_k} & = \frac{\mu}{6Nh}(1+3N-3a)\Delta s_{p_k}\Delta q_{k,a+1}\\&
-\frac{\mu}{3Nh}\Delta s_{p_k}\Delta q_{k,a}
\\&-\frac{\mu}{6Nh}(5+3N-3a)\Delta s_{p_k}\Delta q_{k,a-1}\\&
-\frac{h}{2}\D_{q_{k,a}} V_k -\frac{h}{2} \D_{q_{k,a}} V_{k+1}
,\\
\D_{q_{k,N+1}} L_{d_k} & = -\frac{\mu}{6Nh}\Delta s_{p_k}\Delta q_{k,N+1}
-\frac{\mu}{3Nh}\Delta s_{p_k}\Delta q_{k,N}
\\&-\frac{h}{2}\D_{q_{k,N+1}} V_k -\frac{h}{2} \D_{q_{k,N+1}} V_{k+1}
,\\
\T_I^* \L_{R_k}\cdot \D_{R_k} L_{d_k} & = \frac{M}{h} ((F_k-I)\rho_c)^\wedge R_k^T \Delta q_{k,N+1}
+\frac{h}{2} Mg\hat\rho_c R_k^T e_3\\& + \frac{h}{2} Mg F_k \hat\rho_c F_k^T R_k^T e_3,\\
\D_{s_{p_k}} V_k & = -\mu g r_d\cdot e_3 + \mu g d \sin((s_{p_k}-b)/d)\\&
+\frac{1}{2N}\sum_{a=1}^N  \mu g  e_3\cdot(2r_p + q_{k,a} + q_{k,a+1})\\
&+\frac{EA}{l_k^2} (\norm{q_{k,a+1}-q_{k,a}}^2-l_k^2),\\
\D_{q_{k,a}} V_k & = -\mu g l_k e_3+\nabla V^e_{k,a-1} - \nabla V^e_{k,a},\\
\D_{q_{k,N+1}} V_k & = -(\frac{1}{2}\mu l_k+M)g e_3+\nabla V^e_{k,N}.
\end{align*}

The co-Adjoint map on a real space is the identity map on that real space. The co-Adjoint map on $\SO$ is given by $\Ad^*_{F_k^{-1}} p = F_k p = (F_k \hat p F_k^T)^\vee$ for any $p\in (\Re^3)^*\simeq\so^*$. Using the product structure of the configuration manifold, we have
\begin{align}
\Ad^*_{f_k^{-1}}(\T_e^*\L_{f_{k}}\cdot & \D_{f_{k}}L_{d_{k}}) = \big[ \D_{\Delta s_{p_k}} L_{d_k};\, \D_{\Delta q_{1,k}} L_{d_k},\cdots,\nonumber\\
&\D_{\Delta q_{N+1,k}} L_{d_k};\, \Ad^*_{F_k^T}(\T_I^*\L_{F_{k}}\cdot \D_{F_{k}}L_{d_{k}})\big],\label{eqn:coAdTLLd}
\end{align}
where
\begin{align}
\Ad^*_{F_k^T} (\T_I^* & \L_{F_{k}} \cdot  \D_{F_{k}}L_{d_{k}})=\frac{1}{h} (F_kJ_d -J_dF_k )^\vee\nonumber\\
& + \frac{M}{h} F_k\hat\rho_c F_k^T R_k^T \Delta q_{k,N+1} + \frac{h}{2} MgF_k\hat \rho_c F_k^T R_k^T e_3.\end{align}

\paragraph{Discrete-time Euler-Lagrange Equations}

Substituting the derivatives of the discrete Lagrangian given in \refeqn{DdelqLd}, \refeqn{TLDFLd} and the appendix, the co-Adjoint map given by \refeqn{coAdTLLd}, and the contributions of the external control moment and the Carnot energy loss term \refeqn{Qdk}, \refeqn{udk} into the Lie group variational integrator on an arbitrary Lie group given by \refeqn{DEL_G0}, \refeqn{DEL_G1}, we obtain the discrete-time Euler-Lagrange equations of the string pendulum at \refeqn{DEL0}-\refeqn{DEL4}.

\begin{table*}
{\allowdisplaybreaks
\begin{gather}
\frac{1}{h}M^0_k\Delta s_{p_k}
+\frac{1}{h}\sum_{a=2}^{N} (M^{31}_{k,a}+M^{23}_{k,a-1})\cdot \Delta q_{k,a} 
+\frac{1}{h} M^{23}_{k,N} \Delta q_{k,N+1} 
-\frac{1}{h}M^0_{k-1}\Delta s_{p_{k-1}}
-\frac{1}{h}\sum_{a=2}^{N} (M^{31}_{k-1,a}+M^{23}_{k-1,a-1})\cdot \Delta q_{k-1,a} \nonumber\\
-\frac{1}{h} M^{23}_{k-1,N} \Delta q_{k-1,N+1} 
- \frac{\mu}{3 h}  \Delta s_{p_k}^2 
+ \frac{\mu}{6N h} \sum_{a=1}^N (\Delta q_{k,a}\cdot \Delta q_{k,a}
+\Delta q_{k,a+1}\cdot \Delta q_{k,a+1}
+\Delta q_{k,a}\cdot \Delta q_{k,a+1})+h \D_{s_{p_k}} V_k+\frac{h}{d} u_k\nonumber\\ 
= \frac{h}{2l_k^2} (\mub \Delta s_{p_k}^2 /h^2 +EA) ( \norm{q_{k,2}}-l_k)^2,\label{eqn:DEL0}
\\
q_{k,1}=0,\label{eqn:DEL1}\\
\frac{1}{h} M^{12}_k \Delta q_{k,a-1} 
+ \frac{2}{h} M^1_k \Delta q_{k,a} 
+\frac{1}{h} M^{12}_k \Delta q_{k,a+1} 
+\frac{1}{h}(M^{31}_{k,a}+M^{23}_{k,a-1})\Delta s_{p_k}
-\frac{1}{h} M^{12}_{k-1} \Delta q_{k-1,a-1} 
-\frac{2}{h} M^1_{k-1} \Delta q_{k-1,a} \nonumber
\\
-\frac{1}{h} M^{12}_{k-1} \Delta q_{k-1,a+1} 
-\frac{1}{h}(M^{31}_{{k-1},a}+M^{23}_{k-1,a-1})\Delta s_{p_{k-1}}
-\frac{\mu}{6Nh}(1+3N-3a)\Delta s_{p_k}\Delta q_{k,a+1}
+\frac{\mu}{3Nh}\Delta s_{p_k}\Delta q_{k,a}\nonumber\\
+h\D_{q_{k,a}} V_k
+\frac{\mu}{6Nh}(5+3N-3a)\Delta s_{p_k}\Delta q_{k,a-1}
= 0,\label{eqn:DEL2}
\\
\frac{1}{h} (M^2_k+M) \Delta q_{k,N+1} 
+\frac{1}{h} M_k^{12} \Delta q_{k,N} 
+ \frac{1}{h} M^{23}_{k,N} \Delta s_{p_k}
+\frac{1}{h} MR_k(F_k-I)\rho_c
-\frac{1}{h} (M^2_{k-1}+M) \Delta q_{k-1,N+1} 
-\frac{1}{h} M_{k-1}^{12} \Delta q_{k-1,N} \nonumber\\
- \frac{1}{h} M^{23}_{k-1,N} \Delta s_{p_{k-1}}
-\frac{1}{h} MR_{k-1}(F_{k-1}-I)\rho_c
+\frac{\mu}{6Nh}\Delta s_{p_k}\Delta q_{k,N+1}
+\frac{\mu}{3Nh}\Delta s_{p_k}\Delta q_{k,N}
+h\D_{q_{k,N+1}} V_k = 0,\label{eqn:DEL3}
\\
\frac{1}{h} (F_kJ_d -J_dF_k^T-J_dF_{k-1} +F_{k-1}^TJ_d)^\vee 
+ \frac{M}{h} \hat\rho_c  R_k^T (\Delta q_{k,N+1} -\Delta q_{k-1,N+1})
-h Mg\hat\rho_c R_k^T e_3 =0,\label{eqn:DEL3h}
\\
s_{p_{k+1}} = s_{p_k} +\Delta s_{p_k},\quad q_{k+1,a}=q_{k,a}+\Delta q_{k,a},\quad
R_{k+1}=R_k F_k.\label{eqn:DEL4}
\end{gather}}
\hrule
\end{table*}

Equation \refeqn{DEL2} is satisfied for $a\in\{2,\ldots,N\}$, and \refeqn{DEL4} is satisfied for $a\in\{1,\ldots,N\}$. For given $g_k=(s_{p_k};q_{k,1},\ldots q_{k,N+1};R_k)$, we solve \refeqn{DEL0}-\refeqn{DEL3h} for the relative update $f_k=(\Delta s_{p_k}; \Delta q_{k,1},\ldots$, $\Delta q_{k,N+1};F_k)$. Then, the configuration at the next step $g_{k+1}=(s_{p_{k+1}};q_{k+1,1},\ldots,q_{k+1,N+1};R_{k+1})$ can be obtained by \refeqn{DEL4}. This yields a discrete-time Lagrangian flow map $(g_k,f_k)\rightarrow(g_{k+1},f_{k+1})$, and this is repeated. 

\paragraph{Special Cases} 
If we set $\Delta s_{p_k}\equiv 0$ for any $k$, then the discrete-time Euler-Lagrange equations and Hamilton's equations provide a geometric numerical integrator for a string pendulum model with a fixed unstretched string length, studied in~\cite{LeeLeoPICDC09}. If we chose $\rho_c=0$, then these equations describe the dynamics of an elastic string attached to a point mass and a fixed pivot.

\paragraph{Computational Approach}
These Lie group variational integrators for a string pendulum are implicit: at each time step, we need to solve nonlinear implicit equations to find the relative update $f_k\in\G$. Therefore, it is important to develop an efficient computational approach for these implicit equations.  This computational method should preserve the group structure of $f_k$, in particular, the orthogonal structure of the rotation matrix $F_k\in\SO$. The key idea of the computational approach proposed in this paper is to express the rotation matrix $F_k$ in terms of a vector $c_k\in\Re^3$ using the Cayley transformation~\cite{HaiLub00}:
\begin{align}
F_k = (I+\hat c_k)(I-\hat c_k)^{-1}.\label{eqn:CT}
\end{align}
Since the rotation matrix $F_k$ represents the relative attitude update between two adjacent integration steps, it converges to the identity matrix as the integration step $h$ approaches zero. Therefore, this expression is valid for numerical simulations even though the Cayley transformation is a local diffeomorphism between $\Re^3$ and $\SO$. 

Our computational approach is as follows. The implicit equations for $F_k$ given by \refeqn{DEL3h} are rewritten in terms of a vector $c_k$ using  \refeqn{CT}, and a relative update expressed by a vector $X_k=[\Delta s_{p_k};\Delta q_{k,1},\ldots\Delta q_{k,N+1};c_k]\in \Re\times(\Re^3)^{N+1}\times\Re^3$ is solved by using a Newton iteration. After the vector $X_k$ converges, the rotation matrix $F_k$ is obtained by \refeqn{CT}. 

This computational approach is desirable, since the implicit equations are solved numerically using operations in a linear vector space. The three-dimensional rotation matrix $F_k$ is computed by numerical iterations on $\Re^3$, and its orthogonal structure is automatically preserved by \refeqn{CT}. It has been shown that this computational approach is so numerically efficient that the corresponding computational load is comparable to explicit integrators~\cite{LeeLeoCMDA07}.

\begin{acknowledgements}
This research has been supported in part by National Science Foundation Grants DMS-0714223, DMS-0726263, DMS-0747659, ECS-0244977, CMS-0555797.
\end{acknowledgements}

\bibliography{StrPend}
\bibliographystyle{ieeetran}

\end{document}